\documentclass{LMCS}

\usepackage{enumerate}
\usepackage{hyperref}

\usepackage{amssymb}
\usepackage{amsmath}
\usepackage{latexsym}

\def \dom{{\rm dom}}

\def \graph{{\rm graph}}
\def \id{{\rm id}}
\def \imp {\hspace{5pt}\Longrightarrow\hspace{5pt}}
\def \In{{\subseteq}}
\def \om{{\Sigma^\omega }}
\def \pf{:\hspace{0.6ex}\subseteq \hspace{-0.4ex}}

\def \prf{\sqsubseteq}
\def \range{{\rm range}}
\def \s{{\Sigma^*}}

\def \cal{\mathcal}

\def \Is{{\rm id_*}}

\def\IB{{\mathbb{B}}}
\def\IC{{\mathbb{C}}}
\def\IN{{\mathbb{N}}}
\def\IN{{\mathbb{N}}}
\def\IQ{{\mathbb{Q}}}
\def\IR{{\mathbb{R}}}
\def\IR{{\mathbb{R}}}
\def\IZ{{\mathbb{Z}}}
\newcommand{\an}{\ \wedge\ }

\newcommand{\mto}{\rightrightarrows}

\newcommand{\pproof}{\proof }
\newcommand{\qq}{\qed}

\newcommand{\bb}{ \hspace{-0.76ex}- \hspace{-0.80ex}}

\newcommand{\mmto}{\mbox{
\setlength{\unitlength}{1em}
\begin{picture}(0.4,0)
\makebox(0,0.6){$\mbox{\scriptsize \raisebox{0.083em}{$|$}}   \hspace*{-1.1ex}\mto$}
\end{picture}
}}

\def\doi{7 (2:19) 2011}
\lmcsheading%
{\doi}
{1--21}
{}
{}
{Jan.~\phantom03, 2011}
{Jun.~28, 2011}
{}   

\begin{document}

\title[Computation in analysis]{Turing machines on represented sets, a model of computation for Analysis}


\author[N.~R.~Tavana]{Nazanin R.~Tavana}   
\address{Amirkabir University of Technology, Tehran, Iran} 
\email{nazanin$_-$t@aut.ac.ir}  

\author[K.~Weihrauch]{Klaus Weihrauch} 
\address{University of Hagen, Hagen, Germany}    
\email{Klaus.Weihrauch@FernUni-Hagen.de}  





\begin{abstract}
  \noindent We introduce a new type of generalized Turing machines (GTMs), which are intended as a tool for the mathematician who studies computability in Analysis. In a single tape cell a  GTM can store a symbol, a real number, a continuous real function or a probability measure, for example. The model is based on TTE, the representation approach for computable analysis.
 As a main result we prove that the functions  that are computable via given representations are closed under GTM programming.
This generalizes the well known fact that these functions are closed under composition.
The theorem allows to speak about objects themselves instead of names in algorithms and proofs. By using GTMs for specifying algorithms, many proofs become more rigorous and also simpler and more transparent since the GTM model is very simple and allows to apply well-known techniques from Turing machine theory.
We also show how finite or infinite sequences as names can be replaced by  sets
(generalized representations) on which computability is already defined via representations.
This allows further simplification of proofs.
All of this is done for multi-functions, which are essential in Computable Analysis, and multi-representations, which often allow more elegant formulations.
As a byproduct we show that the computable functions on finite and infinite sequences of symbols are closed under programming with GTMs.
We conclude with examples of application.

\end{abstract}

\maketitle

\section{Introduction}\label{seca}

In 1955 A. Grzegorczyk and D. Lacombe \cite{Grz55,Grz57,Lac55e} proposed a new definition of computable real functions. Their idea became the basis of a general approach to computability in Analysis,  TTE (Type-2 Theory of Effectivity), also called the ``representation approach to computable analysis'' \cite{KW85,Wei00,Sch03,BHW08}.
TTE supplies a uniform method for defining natural computability on a variety of spaces
considered in Analysis such as Euclidean space, spaces of continuous real functions, open, closed or compact subsets of Euclidean space, computable metric spaces, spaces of integrable functions, spaces of probability measures, Sobolev spaces and spaces of distributions. There are various other approaches for studying computability in Analysis \cite[Chapter~9]{Wei00}, but for this purpose, still TTE seems to be the most useful one.

In TTE computability of functions on $\s$, the set of finite words, and $\om$, the set of infinite sequences over a finite alphabet $\Sigma$, is defined explicitly by, for example, ``Type-2 Turing machines''.
Via notations $\nu:\s\to X$ or representations $\delta:\om\to X$,  such ``concrete'' finite or infinite sequences are used as ``names'' for ``abstract'' objects such as real numbers, continuous real functions etc.  A function on the abstract objects is called computable, if it can be realized by a computable function on names.

In ordinary computability theory, for proving computability of a word function $g:(\s)^n\to\s$, in general it is not necessary to write a (usually very long) code of a Turing machine. Instead it suffices to sketch an algorithm that uses some ``simpler'' functions already known to be computable. The method can be formalized by introducing an abstract model of computation, for example, Turing machines (let us call them ``P-machines'') that in addition to the usual statements can use some additional word functions $f:(\s)^n\to \s$ for assignments (``subroutines''). A straightforward proof shows that the computable functions are closed (not only under composition but) under programming with P-machines. More precisely,
the function $f_M$ computed by a P-machine $M$ that uses only computable functions as subroutines is computable, that is, computable by an ordinary Turing machine. Therefore, for proving computability of a function $g$ it suffices to describe
informally a P-machine $M$ that uses only computable functions $f$ as subroutines and to prove that $f_M=g$.

In TTE, the situation is similar. For proving computability of a function on ``abstract'' sets it must be shown that there is a realization that is computable on a Type-2 Turing machine.
Since usually defining or even sketching a concrete Type-2 machine is much too cumbersome, in many articles only algorithms are sketched that use functions on ``abstract'' sets already known to be computable. For a while this method has been applied although its soundness has not been proved.
In \cite{Wei08} the second author has closed this gap by introducing an abstract model of computation for TTE, namely flowcharts with indirect addressing and computable functions on abstract data as subroutines. However, for this model the technical framework of definitions, theorems and proofs has turned out to be complicated and nontransparent such that people preferred to continue with informal arguments rather than applying or mentioning the main results from \cite{Wei08}.

\smallskip
In this article we introduce a very simple model of computation for computable analysis, called here generalized Turing machines. It generalizes the ordinary multi-tape Turing machines with finitely many tapes numbered from $0$ to $L$, finitely many input tapes and one output tape as follows: a generalized Turing machine has a finite tape alphabet $\Gamma$ and for each tape $i$ a set $X_i$. Each cell of Tape $i$ contains an element $x\in \Gamma \cup X_i$. In addition to the usual Turing machine statements on each tape (move left, move right, write $a\in\Gamma$, branch if $a\in\Gamma$ is scanned by the head) two further kinds of statements are allowed  (where $x_j$ is the content of the cell scanned by the head on Tape $i$):

(1) assignments ``$i:=f(i_1,\ldots,i_n)$ '' where $f:X_{i_1}\times \ldots\times X_{i_n}\mto X_i$ is a multi-function (meaning: write some $x\in f(x_{i_1},\ldots x_{i_n})$ on the cell scanned by the head of Tape $i$),

 (2) branchings ``$({\rm if}\ f(i_1,\ldots,i_n) \ {\rm then}\  l', \ {\rm else}\ l'')$'' where
$f\pf X_{i_1}\times\ldots\times X_{i_n}\to \s$ is a partial function (meaning:
if $f(x_{i_1},\ldots x_{i_n})=0$ then go to Label $l'$,
if $f(x_{i_1},\ldots x_{i_n})=1$ then go to Label $l''$, and loop otherwise). ($\Sigma$ will be an alphabet with $0,1\in\Sigma$.)
\smallskip

The model allows to use the universal computational power of Turing machines and for each set $X_j$ used in a machine the number of elements $x\in X_j$ that can be stored during a computation is not bounded. Generalized Turing machines share these properties with the flowcharts with indirect addressing  \cite{Wei08} and with the WhileCC* programs \cite{TZ04}. Our generalized Turing machines can be considered also as a generalization of the BSS-machine \cite{BSS89,BCSS98,Blu04}. In the BSS-model for the real numbers  the algebraic operations and the test ``$x<y$'' are allowed. But in Computable Analysis the test ``$x<y$'' is and should not be computable \cite[Chapter~9]{Wei00}\cite{Bra03f,BC06}.

In Section~\ref{secb} we summarize some mathematical preliminaries, in particular realization of multi-functions by multi-functions via generalized multi-representations.
Generalized multi-representations allow
 simpler but still abstract data as names instead of sequences of
symbols. This generalizes \cite{Bla00} where domains are allowed as sets of names.

The new model of generalized Turing machines and their semantics are defined in Section~\ref{secc}. In Section~\ref{secd} we generalize the concept of multi-representation from sets to machines and prove that realization is not only closed under composition but under programming with generalized Turing machines. In Section~\ref{sece} we prove that for a generalized Turing machine $M$ such that $Y_i\in\{\s,\om\}$ for all $i$ that contains only computable functions on $\s$ and $\om$, the function $f_M$ on $\s$ and $\om$ is computable (accordingly for continuity) (Theorem~\ref{t3} and Corollary~\ref{co2} cf. \cite[Theorem~15]{Wei08}).

The main results are proved in Section~\ref{secf}. If $P$ is a generalized Turing machine where for every tape $i$ the set $Z_i$ is equipped with a multi-representation $\delta_i:\om\mto Z_i$ and every function on the $Z_i$ used in the machine is relatively computable via the corresponding multi-representations, then the function $f_P$ computed by the machine is relatively computable via the corresponding multi-representations
(Theorem~\ref{t4}, cf. \cite[Theorem~30]{Wei08}).
Roughly speaking, the relatively computable functions are closed under programming. The theorem holds accordingly for continuous instead of computable functions.
The theorem holds accordingly if the $\delta_i:Y_i\mto Z_i$ are generalized multi-representations and the functions on the realizing sets $Y_i$ used in the machine are computable w.r.t. a family  $(\gamma_i:\om\mto Y_i)_i$ of multi-representations (Theorem~\ref{t5}, cf. \cite[Theorem~31]{Wei08}).
This theorem allows to use the concept of realization rigorously in a more abstract and often simpler way.
Both theorems allow to formulate and argue about algorithms in terms of ordinary analysis and almost no mentioning of concrete representations.
In Section~\ref{secg} some examples illustrate the main results.
In particular, we present a method for proving the relation $\leq_W$ introduced in
 \cite{BG11} for comparing the non-computability theorems in analysis.
 As an addendum to this introduction the reader is referred to \cite[Section~1]{Wei08}.

\section{Preliminaries}\label{secb}

In this section we summarize some mathematical preliminaries. For more details see \cite{Wei00,Wei08}. Let $\Sigma$ be a non-empty finite set which is called \emph{alpahabet}. We assume $0,1\in\Sigma$. Classically, computability is introduced for functions $f\pf (\s)^n\rightarrow\s$ on the set $\s$ of finite words over $\Sigma$, for example by means of Turing machines. For computing functions on other sets $M$ such as natural numbers, rational numbers and finite graphs, words are used as codes or names of elements of $M$. Under this view a machine transforms words to words without understanding the meaning given to them by the user. We can extend this concept by using infinite sequences of symbols of $\Sigma$ as names and by defining computability for functions which transform such infinite sequences. The set $\om$ of infinite sequences of symbols from $\Sigma$ has the same cardinality as the set of real numbers, therefore it can be used as a set of names for every set with at most continuum cardinality such as real numbers, the set of open subsets of $\mathbb{R}$ and the set $C[0,1]$ of real continuous functions on the interval $[0,1]$.

A \emph{multi-function}  from $A$  to $B$ is a triple
$f = (A,B,R_f )$ such that $R_f\subseteq A\times B$ (the graph of f).
We will denote it by $f : A \mto B$. (The concept of multi-function can be considered as a generalization of the concept of partial function. There is no need for a separate notation  for ``total'' multi-functions.) For $a\in A$,
let $f(a) := \{b\in B | (a, b)\in R_f \}$. For $X\subseteq A$ let
$f[X] := \{b\in B | (\exists a\in X)(a,b)\in R_f \}$, $\dom(f) := \{a\in A | f(a) \neq\emptyset\}$, and $\range(f) := f[A]$.
If, for every $a\in A$, $f(a)$ contains at most one element,
$f$ is a usual partial function denoted by $f\pf A\rightarrow B$. We write ``$f(a)\downarrow$'' ($f(a)$ exists) if $a\in\dom(f)$ and ``$f(a)\uparrow$'' ($f(a)$ diverges) if $a\not\in\dom(f)$ .

In the intended applications, for a multi-function $f: A\mto B$,
$f(a)$ is interpreted as the set of all results which are
``acceptable'' on input $a\in A$. Any concrete computation, a
realization of $f$,  will produce on input $a\in\dom(f)$ some element
$b\in f(a)$, but often there is no method to select a specific one (see \cite{Luc77,BC90,Wei00} and the examples in \cite[Section~3]{Wei08}). The following definition of composition $g\circ f : A\mto C$ of multi-functions $f: A\mto B$ and $g: B\mto C$  is in accordance with this interpretation:
$a\in \dom(g\circ f)$  iff $(a\in\dom(f)\an f(a)\In \dom(g))$
and $g\circ f(a) \; := \; g[f(a)]$ for all $a\in \dom(g\circ f)$.
For the composition of multi-representations we will use the ``relational'' or ``non-deterministic'' composition $\odot$, see (\ref{f23}) in Section~\ref{secf}.

For  $u,v\in \s\cup\om$, $u\prf v$ ($u$ is a {\em prefix} of $v$) iff $v=uw$ for  some $w\in \s\cup\om$. For vectors over $\s\cup\om$ define
$(u_1,\ldots ,u_n)\prf(v_1,\ldots,v_n)$ iff $(\forall i)\,u_i\prf v_i$.
Computable functions on $\s$ can be defined by Turing machines \cite{HU79}.
Computable functions on $\s$ and $\om$ can be defined by Type-2 machines \cite{Wei00}.
A \emph{Type-2 machine} $M$ is a multi-tape Turing machine with $k$ input tapes
(for some $k\geq0$), finitely many work tapes and a single one-way output tape
together with a type specification $(Y_1,\dots,Y_k\rightarrow Y_0)$, $Y_i\in\{\om,\s\}$.

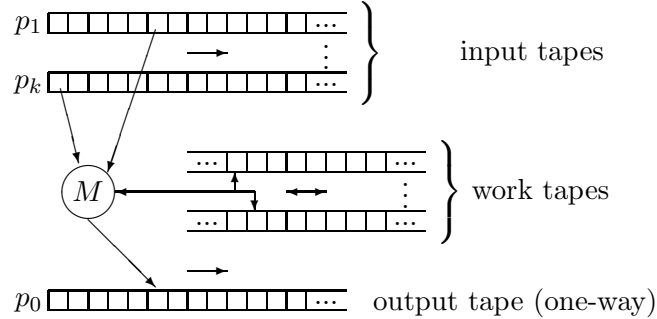
\begin{figure}[htb]
\setlength{\unitlength}{1,5pt}
\linethickness{0.4pt}
\begin{picture}(107.00,83.00)(30,10)
\multiput(10,15)(5,0){14}{\line(0,1){5.00}}
\put(80.00,17.00){\makebox(0,0)[cc]{$...$}}
\put(20.00,45.00){\circle{14.00}}
\put(20.00,45.00){\makebox(0,0)[cc]{$M$}}
\multiput(10,70)(5,0){14}{\line(0,1){5.00}}
\put(80.00,72.00){\makebox(0,0)[cc]{$...$}}
\multiput(10,85)(5,0){14}{\line(0,1){5.00}}
\put(80.00,87.00){\makebox(0,0)[cc]{$...$}}
\multiput(55,35)(5,0){9}{\line(0,1){5.00}}
\put(100.00,37.00){\makebox(0,0)[cc]{$...$}}
\multiput(55,50)(5,0){9}{\line(0,1){5.00}}
\put(100.00,52.00){\makebox(0,0)[cc]{$...$}}
\put(50.00,37.00){\makebox(0,0)[cc]{$...$}}
\put(50.00,52.00){\makebox(0,0)[cc]{$...$}}
\put(45.00,80.00){\vector(1,0){10.00}}
\put(45.00,25.00){\vector(1,0){10.00}}
\put(75.00,45.00){\vector(1,0){5.00}}
\put(75.00,45.00){\vector(-1,0){5.00}}
\put(87.00,80.00){\makebox(0,0)[lc]{$\left.\rule{0mm}{22pt}\right\}\;$ \hspace{4ex}
input tapes}}
\put(107.00,45.00){\makebox(0,0)[lc]{$\left.\rule{0mm}{22pt}\right\}$
 work tapes}}
\put(80.00,81.00){\makebox(0,0)[cc]{$\vdots$}}
\put(100.00,46.00){\makebox(0,0)[cc]{$\vdots$}}
\put(10.00,90.00){\line(1,0){75.00}}
\put(10.00,85.00){\line(1,0){75.00}}
\put(10.00,75.00){\line(1,0){75.00}}
\put(10.00,70.00){\line(1,0){75.00}}
\put(45.00,55.00){\line(1,0){60.00}}
\put(45.00,50.00){\line(1,0){60.00}}
\put(45.00,40.00){\line(1,0){60.00}}
\put(45.00,35.00){\line(1,0){60.00}}
\put(10.00,20.00){\line(1,0){75.00}}
\put(10.00,15.00){\line(1,0){75.00}}
\put(13.00,71.00){\vector(1,-4){4.67}}
\put(37.00,86.00){\vector(-1,-3){12.00}}
\put(20.00,38.00){\vector(1,-1){17.00}}
\put(62.00,45.00){\vector(-1,0){35.00}}
\put(62.00,45.00){\vector(0,-1){5.00}}
\put(57.00,45.00){\vector(0,1){5.00}}
\put(5.00,87.00){\makebox(0,0)[cc]{$p_1 $}}
\put(5.00,72.00){\makebox(0,0)[cc]{$ p_k$}}
\put(5.00,17.00){\makebox(0,0)[cc]{$p_0$}}
\put(87.00,17.00){\makebox(0,0)[lc]{ \ output tape (one-way)}}
\end{picture}
\caption{A Type-2 machine} \label{fig2}
\end{figure}
\noindent The function $f_M \pf Y_1\times\dots \times Y_k\rightarrow Y_0$ computed by the Type-2 machine $M$ is
defined as follows:\\
Case $Y_0=\s$: $f_M(p_1, \dots, p_k) = w$, iff $M$ halts on input $(p_1,\dots , p_k)$ with
$w\in\s$ on the output tape;\\
Case $Y_0=\om$: $f_M(p_1,\dots, p_k) = p_0$, iff $M$ computes forever on input
$(p_1,\dots , p_k)$ and writes $p_0\in\om$ on the output tape.\\
We call a function $f\pf Y_1\times\dots \times Y_k\rightarrow Y_0$ Turing computable, iff $f = f_M$ for some Type-2 machine $M$, and deviant from the usual terminology
we call it computable, if it has a Turing computable extension. (Notice that usually
``computable'' means Turing computable.) The computable functions are closed under composition.

On $\s$ we consider the discrete topology and on $\om$ the Cantor topology defined by the basis $\{u\om\mid u\in\s\}$ of open sets. As a fundamental result, every computable function on $\s$ and $\om$ is continuous.

A {\em representation} of a set $M$ is a surjective function $\delta\pf Y\to M$ where
$Y=\s$ or $Y=\om$. (Often the word ``representation'' is reserved for the case $\delta\pf\om\to M$ and surjective functions $\nu\pf\s\to M$ are called ``notations'' \cite{Wei00}.) We will use {\em multi-representations} $\delta: Y\mto M$ where $M=\range(\delta)$ ($Y\in\{\s,\om\}$). Here, a name $w\in y$ may be a name of many $x\in M$. Finally we use \emph{generalized multi-representations} $\lambda: U\mto M$ such that $\range(\lambda)=M$ where an arbitrary set $U$ is considered as the set of ``names''.

If for a generalized multi\bb representation
$\delta:U\mto X$,  $x\in\delta(u)$ then we say ``$u$ realizes $x$ (via $\delta$)'' or
``$u$ is a name of $x$''. The {\em realization} of functions by functions is a central concept in TTE. We define the most general case: the realization of a multi-function by  a multi-function via generalized multi-representations.

\begin{defi}[realization]\cite{Wei08}\label{d3}
Let $f:X_1\times\ldots\times X_n\mto X_0$ and $g:  Y_1\times \ldots\times Y_n\mto Y_0$
be multi-functions and let $\gamma_i:X_i\mto Y_i$ ($0\leq i \leq n$) be generalized multi-representations.
For $x=(x_1,\ldots,x_n)\in X_1\times\ldots\times X_n$ let $\gamma(x):=\gamma_1(x_1)\times \ldots\times \gamma_n(x_n)$.

Then ``$f$ realizes  $g$ via $(\gamma_1,\ldots,\gamma_n,\gamma_0)$'' or ``$f$ is a $(\gamma_1,\ldots,\gamma_n,\gamma_0)$-realization of $g$'', iff
for all $x\in X_1\times\ldots\times X_n$ and
$y\in Y_1\times \ldots\times Y_n$,
\begin{eqnarray}\label{f1}
y\in \gamma(x)\cap \dom(g) \  \Longrightarrow \  \big(\,f(x)\neq\emptyset\an
(\forall x_0\in f(x))\, g(y)\cap \gamma_0(x_0) \not=\emptyset\,\big)\,.
\end{eqnarray}
\end{defi}

 \begin{figure}[ht]
\unitlength0.26ex
\begin{picture}(110,70)(0,-8)

\put(0,0){\circle*2}
\put(0,50){\circle*2}
\put(70,0){\circle*2}
\put(70,50){\circle*2}

\put(0,0){\vector(1,0){69}}
\put(0,0){\vector(4,1){20}}
\put(0,0){\vector(4,-1){25}}

\put(0,50){\vector(0,-1){49}}
\put(0,50){\vector(-1,-4){5}}
\put(0,50){\vector(1,-4){4}}

\put(0,50){\vector(1,0){69}}
\put(0,50){\vector(4,1){22}}
\put(0,50){\vector(4,-1){24}}

\put(70,50){\vector(0,-1){49}}
\put(70,50){\vector(-1,-4){5}}
\put(70,50){\vector(1,-4){4.5}}

\put(-7,50){\makebox(5,8){\boldmath$x$}}
\put(70,50){\makebox(14,8){\boldmath$x_0$}}
\put(-7,-8){\makebox(5,8){\boldmath$y$}}
\put(80,-10){\makebox(46,8) {\boldmath $\,y_0\in g(y)\cap\gamma_0(x_0)$}}

\put(33,52){\makebox(5,8){\boldmath$f$}}
\put(33,-10){\makebox(5,8){\boldmath$g$}}
\put(-13,21){\makebox(5,8){\boldmath$\gamma$}}
\put(77,21){\makebox(5,8){\boldmath$\gamma_0$}}

\end{picture}
\caption { $f$ realizes $g$ via $(\gamma_1,\ldots,\gamma_n,\gamma_0)$.}
\label{bildc}
\end{figure}
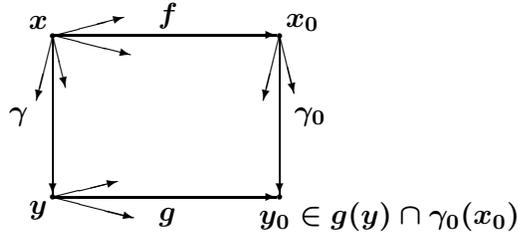
\noindent Figure~\ref{bildc} illustrates the realization of $g$ by
$f$. Roughly speaking, provided $x$ is a $\gamma$-name of $y\in\dom(g)$ then
\[\begin{array}{rlrllllll}
f(x) & \mbox{is a name of} & g(y)&  \mbox{if} & f & \mbox{is single-v. and } &  g & \mbox{is single-v.}\, ,\\
f(x) & \mbox{is a name of some} & y_0\in g(y)&  \mbox{if} &f & \mbox{is single-v. and } &  g & \mbox{is multi-v.}\, ,\\

\mbox{every}\ \ x_0\in f(x)
 & \mbox{is a name of some} & y_0\in g(y)&  \mbox{if} & f & \mbox{is multi-v. and } &  g & \mbox{is multi-v.}\, .
\end{array}
\]
For further technical details see \cite {Wei00} and \cite[Sections 1,2,3,6,8 (until Lemma~28) and 9]{Wei08}.

\section{Generalized Turing machines}\label{secc}
We generalize multi-tape Turing machines \cite{HU79} to generalized Turing machines as follows.
A  generalized Turing machine ({\rm GTM}) has $L+1$ tapes where Tapes $1,\ldots, k$ are the input tapes, Tapes $k+1,\ldots, L$ are work tapes and
Tape $0$ is the output tape. There is a finite work alphabet $\Gamma$ and the blank symbol $b\in\Gamma$. For an ordinary Turing machine, there is a finite input/output alphabet $\Sigma$ such that $\Sigma\cap \Gamma=\emptyset$ and at any time every cell of every tape contains exactly one element (``symbol'') $a\in\Sigma\cup \Gamma$. We generalize the definition by assigning to every tape $i$ a set $X_i$ (which may be empty) such that at any time every cell of Tape $i$ contains exactly one element  $a\in X_i\cup \Gamma$.

\begin{figure}[h]

\setlength{\unitlength}{1.8pt}
\linethickness{0.7pt}
\begin{picture}(110,90)(-00,-5)

\multiput(47.5,80)(0,-2){44}{\line(0,-1){1.2}}
\multiput(54.5,80)(0,-2){44}{\line(0,-1){1.2}}

\put(15,76){\makebox(0,0){\footnotesize $-6$}}
\put(21,76){\makebox(0,0){\footnotesize $-5$}}
\put(27,76){\makebox(0,0){\footnotesize $-4$}}
\put(33,76){\makebox(0,0){\footnotesize $-3$}}
\put(39,76){\makebox(0,0){\footnotesize $-2$}}
\put(45,76){\makebox(0,0){\footnotesize $-1$}}
\put(51.5,76){\makebox(0,0){\footnotesize $0$}}
\put(57.5,76){\makebox(0,0){\footnotesize $1$}}
\put(63.5,76){\makebox(0,0){\footnotesize $2$}}
\put(69.5,76){\makebox(0,0){\footnotesize $3$}}
\put(75.5,76){\makebox(0,0){\footnotesize $4$}}
\put(81.5,76){\makebox(0,0){\footnotesize $5$}}
\put(87.5,76){\makebox(0,0){\footnotesize $6$}}

\put(115,76){\sf\small cell number}
\put(-33,76){\parbox{33pt}{\sf\small \hspace{1.5ex}tape number}}
\put(-25,64){\makebox(0,0)[cc]{$0$}}
\put(-25,46){\makebox(0,0)[cc]{$1$}}
\put(-25,30){\makebox(0,0)[cc]{$k$}}
\put(-25,18){\makebox(0,0)[cc]{$k+1$}}
\put(-25,2){\makebox(0,0)[cc]{$L$}}
\put(-25,40){\makebox(0,0)[cc]{$ \vdots$}}
\put(-25,12){\makebox(0,0)[cc]{$ \vdots$}}

\put(-16,76){\parbox{33pt}{\sf\small \hspace{1.5ex}set}}
\put(-10,64){\makebox(0,0)[cc]{$X_0$}}

\put(-10,46){\makebox(0,0)[cc]{$X_1$}}
\put(-10,30){\makebox(0,0)[cc]{$X_k$}}
\put(-10,18){\makebox(0,0)[cc]{$X_{k+1}$}}
\put(-10,2){\makebox(0,0)[cc]{$X_L$}}
\put(-10,40){\makebox(0,0)[cc]{$ \vdots$}}
\put(-10,12){\makebox(0,0)[cc]{$ \vdots$}}

\put(0,62){\line(1,0){102}}
\put(0,68){\line(1,0){102}}
\multiput(12,62)(6,0){14}{\line(0,1){6}}
\put(6,64.5){\makebox(0,0)[cc]{$.\,.\,.$}}
\put(96,64.5){\makebox(0,0)[cc]{$.\,.\,.$}}
\put(115,63.5){\parbox{27pt}{\sf\small output  tape}}

\put(0,44){\line(1,0){102}}
\put(0,50){\line(1,0){102}}
\multiput(12,44)(6,0){14}{\line(0,1){6}}
\put(6,46.5){\makebox(0,0)[cc]{$.\,.\,.$}}
\put(96,46.5){\makebox(0,0)[cc]{$.\,.\,.$}}

\put(63,40){\makebox(0,0)[cc]{ $\Huge \vdots$}}

\put(0,28){\line(1,0){102}}
\put(0,34){\line(1,0){102}}
\multiput(12,28)(6,0){14}{\line(0,1){6}}
\put(6,30.5){\makebox(0,0)[cc]{$.\,.\,.$}}
\put(96,30.5){\makebox(0,0)[cc]{$.\,.\,.$}}
\put(105,38.7){\makebox(0,0)[lc]{$\left.\rule{0pt}{25pt}\right\}$ }}
\put(115,38){\parbox{27pt}{\sf\small input  tapes}}

\put(0,16){\line(1,0){102}}
\put(0,22){\line(1,0){102}}
\multiput(12,16)(6,0){14}{\line(0,1){6}}
\put(6,18.5){\makebox(0,0)[cc]{$.\,.\,.$}}
\put(96,18.5){\makebox(0,0)[cc]{$.\,.\,.$}}

\put(63,12){\makebox(0,0)[cc]{$\Huge\vdots$}}

\put(0,0){\line(1,0){102}}
\put(0,6){\line(1,0){102}}
\multiput(12,0)(6,0){14}{\line(0,1){6}}
\put(6,2.5){\makebox(0,0)[cc]{$.\,.\,.$}}
\put(96,2.5){\makebox(0,0)[cc]{$.\,.\,.$}}
\put(105,10.7){\makebox(0,0)[lc]{$\left.\rule{0pt}{25pt}\right\}$}}
\put(115,10){\parbox{38pt}{\sf\small additional work tapes}}

\newsavebox{\head}
\savebox{\head}{\linethickness{1.3pt}
\put(0,0){\line(1,0){8}}
\put(0,0){\line(0,1){8}}
\put(0,8){\line(1,0){8}}
\put(8,0){\line(0,1){8}}
}

\put(47,61){\usebox{\head}}
\put(35,43){\usebox{\head}}
\put(83,27){\usebox{\head}}
\put(17,15){\usebox{\head}}
\put(65,-1){\usebox{\head}}

\end{picture}
\caption{A generalized Turing machine.} \label{fig1}
\end{figure}
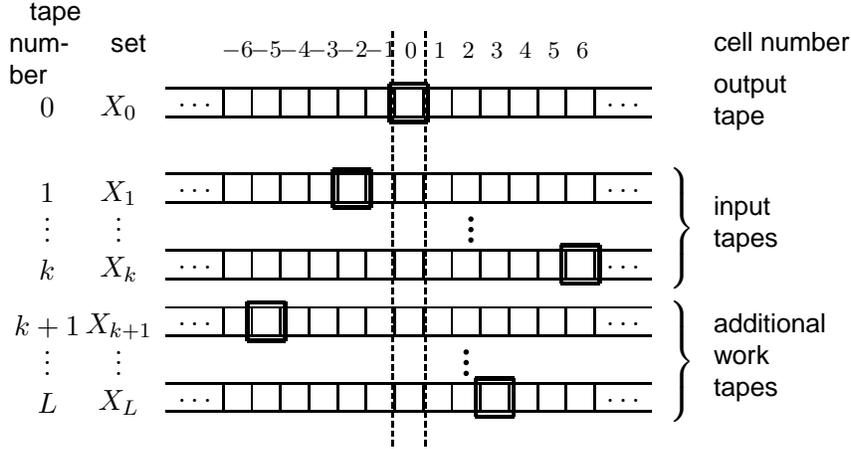

As for an ordinary Turing machine every tape has a read/write head that scans exactly one cell and  there is a finite set ${\mathcal L}$ of labels (usually called states) with an initial label $l_0\in {\mathcal L}$ and a final label $l_f\in {\mathcal L}$.
For every label $l\neq l_f$ there is a statement defining some action on some tape and the next label. As for an ordinary Turing machine in one step on some tape the head can be moved one position to the right or to the left, and for every  symbol $a\in\Gamma$, $a$ can be written on the cell scanned by the head and it can be tested whether  $a$ is scanned by the head (branching). Figure~\ref{fig1} shows the tapes and heads of a generalized Turing machine.

\noindent
Generalized Turing machines may have a further kind of assignments and a further kind of branchings. Let~$x_i$ be the content of the cell scanned by the head on Tape $i$ ($0\leq i\leq L$).
\begin{enumerate}[(1)]
\item \label{e1} ``$(i:=f(i_1,\ldots, i_n),l')\,$'' for some $f:X_{i_1}\times\ldots\times X_{i_n}\mto X_{i}$ meaning: \\
write some $ y\in f(x_{i_1},\ldots x_{i_n})$ on the cell scanned by the head on Tape $i$ and then go to Label $l'$;

\item \label{e3} ``$({\rm if}\ f(i_1,\ldots,i_n) \ {\rm then}\  l', \ {\rm else}\ l'')$''
for some $f\pf X_{i_1}\times\ldots\times X_{i_n}\to \s$ meaning:\\
if $f(x_{i_1},\ldots x_{i_n})=0$ then go to Label $l'$,
if $f(x_{i_1},\ldots x_{i_n})=1$ then go to Label $l''$ (and loop otherwise).
\end{enumerate}

\begin{defi}\label{d1}
A generalized Turing machine (\,{\rm GTM}) is  a tuple \\
${\bf M}=({\mathcal L},l_0,l_f,\Gamma, b, k,L,(X_i)_{0\leq i\leq L},{\rm Stm})$ such that:
\begin{enumerate}[(1)]
\item \label{d1a} ${\cal L}$ is a finite set (``labels''), $l_0,l_f\in {\cal L}$ (``initial'' and ``final'' label);
\item \label{d1b} $\Gamma$ (``work alphabet'') is a finite set,  $\Sigma\cap\Gamma=\emptyset$ and $b\in\Gamma$ (``blank'' symbol);
\item \label{d1c} $k,L\in\IN$, $k \leq L$ (\,$0,1,\ldots, L$: numbers of the tapes; $1,\ldots, k$: numbers of the input tapes; $0$: number of the output tape);

\item \label{d1d}  $X_i$ is a set such that $X_i\cap \Gamma=\emptyset$ ($\;0\leq i\leq L$);

\item \label{d1e} ${\rm Stm}$ is a function assigning to every label $l\in {\cal L } \setminus\{l_f\}$ a statement from the following list  (\,where $\{i,i_1,\ldots,i_n\}\In\{0,1,\ldots,L\}$ and $l',l''\in {\cal L}$):
\begin {enumerate}[(a)]
\item \label{d1e1} $(i,{\rm right},l')$,
\item \label{d1e2} $(i,{\rm left},l')$,
\item \label{d1e3} $(i:=a, l')$ \ (for some $a\in\Gamma$),
\item \label{d1e4} $(i, {\rm if}\ a\ {\rm then} \ l',\ {\rm else}\ l'')$ (for some $a\in\Gamma$),
\item \label{d1e5} $(i:=f(i_1,\ldots, i_n), l')$ (for some $f:X_{i_1}\times\ldots\times X_{i_n}\mto X_{i}$);

\item \label{d1e7} $({\rm if}\ f(i_1,\ldots,i_n) \ {\rm then}\  l', \ {\rm else}\ l'')$ (for some $f\pf X_{i_1}\times\ldots\times X_{i_n}\to \s$ with $\range(f)\In \{0,1\}$).

\end{enumerate}
\end{enumerate}
\end{defi}
\noindent
Notice that for assignments (\ref{d1e5}) we allow multi-valued functions while for tests (\ref{d1e7}) the functions must be single-valued but may still be partial. For defining the semantics we formalize the tape $i$ with inscription by a function $\alpha_i:\IZ\to X_i\cup\Gamma$ and the head position by a number $m_i\in\IZ$.
In the branching (\ref{d1e7}) we will interpret $0\in\s$ as true and $1\in\s$ as false.

\begin{defi}[semantics]\label{d2} Let ${\bf M}=({\cal L},l_0,l_f,\Gamma, b, k,L,(X_i)_{0\leq i\leq L},{\rm Stm})$ be a generalized Turing machine.
\begin{enumerate}[(1)]
\item \label{d2a} Let ${\cal S}:=\prod_{i=0}^L((X_i\cup \Gamma)^\IZ \times \IZ)$ be the set of {\em states} and ${\cal K}:={\cal L}\times {\cal S}$ be the set of {\em configurations}.
For a configuration $\kappa=(l,(\alpha_0,m_0),\ldots ,(\alpha_L,m_L))$ define local modifications of $\kappa$ as follows:

\begin{tabular}{ll}
$\kappa[{\rm label}\leftarrow l_1]$:  &in $\kappa$ replace the label by $l_1$\\
$\kappa[{\rm head}_i\leftarrow m]$: &in $\kappa$ move the head on Tape $i$ to Position $m$,\\
$\kappa[{\rm cell}_i\leftarrow x]$: &in $\kappa$ write $x$ under the head of Tape $i$ .
\end{tabular}

\item \label{d2c} We define a successor relation $\vdash\In {\cal K\times \cal K}$. Let $\kappa=(l,(\alpha_0,m_0),\ldots ,(\alpha_L,m_L))$ and  $x_j:=\alpha_j(m_j)$ for $0\leq j\leq L$. The successors of $\kappa$ are determined by the statement ${\rm Stm}(l)$ as follows: $\kappa\vdash \kappa'$ iff:

\begin{enumerate}[(a)]
\item \label{d2c1} {\boldmath ${\rm Stm}(l)=(i,{\rm right},l')$:} \  $\kappa'=\kappa[{\rm head}_i\leftarrow m_i+1]\,[{\rm label}\leftarrow l']$,

\item  \label{d2c2}{\boldmath ${\rm Stm}(l)=(i,{\rm left},l')$:} \  $\kappa'=\kappa[{\rm head}_i\leftarrow m_i-1]\,[{\rm label}\leftarrow l']$,

\item \label{d2c3} {\boldmath ${\rm Stm}(l)=(i:= a, l')$:} \
$\kappa'=\kappa[{\rm cell}_i\leftarrow a] \,[{\rm label}\leftarrow l']$  \ \ ,

\item \label{d2c4} {\boldmath ${\rm Stm}(l)=(i, {\rm if}\ a\ {\rm then} \ l',\ {\rm else}\ l'')$:} \\
 $\kappa'=\kappa[{\rm label}\leftarrow l']$ if $x_i=a$, and
$\kappa'=\kappa[{\rm label}\leftarrow l'']$ if $x_i\neq a$,

\item \label{d2c5} {\boldmath ${\rm Stm}(l)=(i:=f(i_1,\ldots, i_n), l')$:}  \\
$\kappa'=\kappa[ {\rm cell}_i\leftarrow x]\,[{\rm label}\leftarrow l']$ for some
$x\in f(x_{i_1},\ldots,x_{i_n})$,

\item \label{d2c7} {\boldmath ${\rm Stm}(l)=({\rm if}\ f(i_1,\ldots,i_n) \ {\rm then}\  l', \ {\rm else}\ l'')$:} \\
$\kappa'=\kappa[{\rm label}\leftarrow l']$ if $f(x_{i_1},\ldots,x_{i_n})=0$, and \\
$\kappa'=\kappa[{\rm label}\leftarrow l'']$ if $f(x_{i_1},\ldots,x_{i_n})=1$.
\end{enumerate}

\item \label{d2d}
A {\em computation} is a (finite or infinite) sequence $(\kappa^0,\kappa^1,\ldots)$ of configurations such that $\kappa^i\vdash \kappa^{i+1}$. A computation is {\em maximal} if it is infinite or its last configuration has no $\vdash$-successor.
A configuration $\kappa=(l,(\alpha_0,m_0),\ldots ,(\alpha_L,m_L))$ is {\em accepting} if $l=l_f$ and $\alpha_0(0)\in X_0$.
An {\em accepting computation} is a finite computation
$(\kappa_0,\kappa_1,\ldots ,\kappa_n)$ such that $\kappa_n$ is accepting.

\item \label{d2e}
For $(x_1,\ldots,x_k)\in X_1\times\ldots \times X_k$ define the initial configuration by
\[{\rm IC}(x_1,\ldots,x_k):=(l^0,(\alpha^0_0,0), \ldots, (\alpha^0_L,0))\]
where $l^0=l_0$,
$\alpha^0_i(0)=x_i$ for $1\leq i\leq k$ and $\alpha^0_i(j)=b$ for all other $(i,j)$.
For every configuration $\kappa=(l,(\alpha_0,m_0),\ldots ,(\alpha_L,m_L))$ define
\[{\rm OC}(\kappa):=\left\{
\begin{array}{ll}
\alpha_0(0) & \mbox{if}\ \ \alpha_0(0)\in X_0\\
{\rm div}  & \mbox{otherwise}.
\end{array}\right.\]
Define the multi-function $f_M:X_1\times\ldots\times X_k\mto X_0$ computed by $\bf M$ as follows: For $x_i\in X_i$ ($\,0\leq i\leq k$) let $x_0\in f_M(x_1,\ldots,x_k)$ iff  (\ref{d2e1}) and (\ref{d2e2}):
\begin{enumerate}[(a)]
\item \label{d2e1} every maximal computation with first configuration ${\rm IC}(x_1,\ldots,x_k)$ is accepting,

\item \label{d2e2}
there exists an accepting computation $(\kappa^0,\ldots,\kappa^n)$ with first configuration $\kappa^0={\rm IC}(x_1,\ldots,x_k)$ such that $x_0={\rm OC}(\kappa^n)$.
\end{enumerate}
\end{enumerate}
\end{defi}

\noindent For input $(x_1,\ldots,x_k)\in X_1\times\ldots \times X_k$, the initial configuration has the label $l_0$, on every tape the head is on position $0$, on the input tape $i$
($1\leq i\leq k$) the cell $0$ contains the value $x_i$, and all other tape cells contain the blank symbol $b\in\Gamma$. In every assignment step (\ref{d2c5}) every
$x\in f(x_{i_1},\ldots,x_{i_n})$ can be chosen.
The result of an accepting computation is the inscription of the cell $0$ on Tape $0$, which must be in $X_0$. A value $x\in X_0$ is in $f(x_1,\ldots,x_n)$, if there is an accepting computation with result $x$ and every maximal computation on the same input is accepting.

\section{Realization is Closed Under Programming}\label{secd}

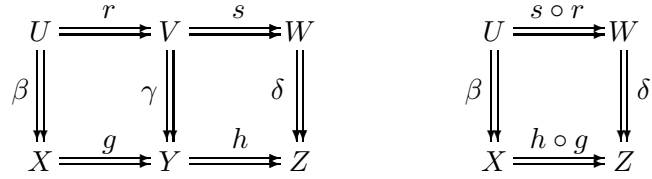
\begin{figure}[htbp]
\unitlength0.26ex
\begin{picture}(40,50)(70,2)

\put(0,40){\makebox(2,2){$U$}}
\put(40,40){\makebox(2,2){$V$}}
\put(80,40){\makebox(2,2){$W$}}

\put(0,0) {\makebox(2,2){$X$}}
\put(40,0){\makebox(2,2){$Y$}}
\put(80,0){\makebox(2,2){$Z$}}

\put(7,0){\vector(1,0){28}}
\put(7,40){\vector(1,0){28}}

\put(7,2){\vector(1,0){28}}
\put(7,42){\vector(1,0){28}}

\put(0,35){\vector(0,-1){28}}
\put(40,35){\vector(0,-1){28}}

\put(2,35){\vector(0,-1){28}}
\put(42,35){\vector(0,-1){28}}

\put(47,0){\vector(1,0){28}}
\put(47,40){\vector(1,0){28}}

\put(47,2){\vector(1,0){28}}
\put(47,42){\vector(1,0){28}}

\put(80,35){\vector(0,-1){28}}
\put(82,35){\vector(0,-1){28}}

\put(20,45){$r$}
\put(20,5){$g$}
\put(60,45){$s$}
\put(60,5){$h$}

\put(-8,20){$\beta$}
\put(32,20){$\gamma$}
\put(72,20){$\delta$}

\put(140,40){\makebox(2,2){$U$}}
\put(180,40){\makebox(2,2){$W$}}
\put(140,0) {\makebox(2,2){$X$}}
\put(180,0){\makebox(2,2){$Z$}}

\put(147,0){\vector(1,0){28}}
\put(147,40){\vector(1,0){28}}
\put(147,2){\vector(1,0){28}}
\put(147,42){\vector(1,0){28}}

\put(140,35){\vector(0,-1){28}}
\put(180,35){\vector(0,-1){28}}
\put(142,35){\vector(0,-1){28}}
\put(182,35){\vector(0,-1){28}}

\put(47,0){\vector(1,0){28}}
\put(47,40){\vector(1,0){28}}

\put(47,2){\vector(1,0){28}}
\put(47,42){\vector(1,0){28}}

\put(80,35){\vector(0,-1){28}}
\put(82,35){\vector(0,-1){28}}

\put(152,45){$s\circ r$}
\put(152,5){$h\circ g$}

\put(132,20){$\beta$}
\put(185,20){$\delta$}
\end{picture}

\caption {$s\circ r$ realizes $h\circ g$, if $r$ realizes $g$ and $s$ realizes $h$.}
\label{bildi}
\end{figure}

For multi-functions on multi-represented sets realization is closed under composition, that is, the composition of realizations realizes the composition \cite[Theorem~3.1.6]{Wei00} \cite[Lemma~20]{Wei08}
 (see Figure~\ref{bildi}).

\noindent Theorem~\ref{t1} generalizes this fact from simple composition to generalized Turing machines. It is the GTM-version of \cite[Theorem~23]{Wei08}.
Let $\Is:\s\to\s$ be the identity on $\s$.
Then by (\ref{f1}), a branching
$f\pf X_1\times\ldots\times X_n\to \s$ is a $(\gamma_1,\ldots,\gamma_n,\Is)$-realization of a branching $g\pf Y_1\times\ldots\times Y_n\to \s$ iff
\begin{eqnarray}\label{f6}
y\in \gamma(x)\cap \dom(g) \  \Longrightarrow \ f(x)= g(y)\,.
\end{eqnarray}
We generalize the concept of realization (Definition~\ref{d3}) from functions to generalized Turing machines as follows:

\begin{defi}\label{d5}
Let ${\bf M}=({\cal L},l_0,l_f,\Gamma, b, k,L,(X_i)_{0\leq i\leq L},{\rm Stm}_M)$ and\\
${\bf N}=({\cal L},l_0,l_f,\Gamma, b, k,L,(Y_i)_{0\leq i\leq L},{\rm Stm}_N)$
be generalized Turing machines and let $\gamma_i:X_i\mto Y_i$ ($\;0\leq i\leq L$) be generalized multi-representations.

Then ``$\bf M$ is a $(\gamma_i)_{i=0}^L$-realization of $\bf N$'' or ``$\bf M$ realizes $\bf N$ via $(\gamma_i)_{i=0}^L$'', \\
if  (\ref{d5a}) -- (\ref{d5d}) for all labels $l\in\cal L$.
\begin{enumerate}[(1)]
\item \label{d5a}
if $\;{\rm Stm}_M(l)\in\{ (i,{\rm right},l'),\ (i,{\rm left},l'),\ (i:= a,l'),\ (i, {\rm if}\ a\ {\rm then} \ l',\ {\rm else}\ l'')\}$ then
$\;{\rm Stm}_M(l) ={\rm Stm}_N(l)$,

\item \label{d5b} if $\;{\rm Stm}_M(l)=(i:=f(i_1,\ldots, i_n), l')$ then
$\;{\rm Stm}_N(l)=(i:=g(i_1,\ldots, i_n), l')$ such that
$f:X_{i_1}\times\ldots\times X_{i_n}\mto X_{i}$ is a
$(\gamma_{i_1},\ldots, \gamma_{i_n},\gamma_i)$-realization of
$g:Y_{i_1}\times\ldots\times Y_{i_n}\mto Y_{i}$,

\item \label{d5d}
if
$\;{\rm Stm}_M(l)=({\rm if}\ f(i_1,\ldots,i_n) \ {\rm then}\  l', \ {\rm else}\ l'')$
 then \\
$\;{\rm Stm}_N(l)=({\rm if}\ g(i_1,\ldots,i_n) \ {\rm then}\  l', \ {\rm else}\ l'')$
such that \\
$f\pf X_{i_1}\times\ldots\times X_{i_n}\to\s$ is a $(\gamma_{i_1},\ldots, \gamma_{i_n}, \Is)$-realization of $g\pf Y_1\times\ldots\times Y_k\to\s$.
\end{enumerate}
\end{defi}

\begin{thm}\label{t1}
Let ${\bf M}=({\cal L},l_0,l_f,\Gamma, b, k,L,(X_i)_{0\leq i\leq L},{\rm Stm}_M)$ and\\
${\bf N}=({\cal L},l_0,l_f,\Gamma, b, k,L,(Y_i)_{0\leq i\leq L},{\rm Stm}_N)$
be generalized Turing machines and let $\gamma_i:X_i\mto Y_i$ ($\;0\leq i\leq L$) be generalized multi-representations.

If $\,\bf M$ realizes $\bf N$ via $(\gamma_i)_{i=0}^L$, then
$f_M:X_1\times\ldots\times X_k\mto X_0$ realizes $f_N:Y_1\times\ldots\times Y_k\mto Y_0$ via $(\gamma_1,\ldots, \gamma_k,\gamma_0)$.
\end{thm}

First we prove a lemma that considers all the details of the generalized Turing machines. It extends the concept of realization for multi-functions in Definition~\ref{d3} to the successor relations $\vdash_M$ and $\vdash_N$.
For a configuration
$\kappa=(l,(\alpha_0,m_0),\ldots ,(\alpha_L,m_L))$ of $\bf M$ and a configuration $\lambda=
 (\overline l,(\overline\alpha_0,\overline m_0),\ldots ,(\overline \alpha_L,\overline m_L)) $ of $\bf N$ we say ``$\kappa$ realizes $\lambda$'' iff (\ref{f4}) and  (\ref{f5}) are satisfied:
\begin{eqnarray}
\label{f4}& l=\overline l \an  (\forall i\in \{0,\ldots, L\})\,m_i=\overline m_i\,, \\  \label{f5} &(\forall \ i\in \{0,\ldots, L\})\,(\forall \  j\in \IZ)\: (\alpha_i(j)=\overline \alpha_i(j)\in\Gamma \,\vee \overline \alpha_i(j)\in\gamma_i\circ \alpha_i(j))\,.
\end{eqnarray}

\begin{lem}\label{l1} Let $\bf M$ be a $(\gamma_i)_{i=0}^L$-realization of $\bf N$. If $\kappa$ realizes $\lambda$ and $\lambda$ has a $\vdash_N$-successor $\lambda''$, then
\begin{enumerate}[\em(1)]
\item \label{l1a} $\kappa$ has a $\vdash_M$-successor $\kappa''$ and
\item \label{l1b} if $\kappa\vdash_M\kappa'$ then there is some $\lambda'$ such that $\lambda\vdash _N\lambda'$  and $\kappa'$ realizes $\lambda'$.
\end{enumerate}
\end{lem}

\pproof Suppose $\kappa$ realizes $\lambda$. By (\ref{f4}) and  (\ref{f5}) $\kappa$ and $\lambda$ can be written as
\[\begin{array}{rcl}
\kappa  & = & (l,(\alpha_0,m_0),\ldots,(\alpha_L,m_L))\\
\lambda & = & (l,(\overline\alpha_0,m_0),\ldots,(\overline\alpha_L,m_L))
\end{array}\]
such that $\alpha_i(j)=\overline \alpha_i(j)\in\Gamma$ or $\overline \alpha_i(j)\in\gamma_i\circ \alpha_i(j))$ for all $i,j$. By assumption, $\lambda$ has a successor $\lambda''$. Then $l\neq l_f$.
We study successively the 6 cases for ${\rm Stm}_N(l)$ from Definition~\ref{d2}.\ref{d2c}.
In the last two cases below let
\[\widetilde x:=(\alpha_{i_1}(m_{i_1}),\ldots,\alpha_{i_n}(m_{i_n}))
  \quad\hbox{and}\quad
\widetilde
y:=(\overline\alpha_{i_1}(m_{i_1}),\ldots,\overline\alpha_{i_n}(m_{i_n})).
\]

{\boldmath ${\rm Stm}_M(l)=(i,{\rm right},l')$:}\\
Then ${\rm Stm}_N(l)=(i,{\rm
right},l')$. By Definition~\ref{d2}.\ref{d2c1} $\lambda''=\lambda[{\rm
head}_i\leftarrow m_i+1][{\rm label}\leftarrow l']$.
Let $\kappa'':=\kappa[{\rm head}_i\leftarrow m_i+1][{\rm
label}\leftarrow l']$.
Then $\kappa\vdash_M \kappa''$ and $\kappa''$ realizes  $\lambda''$.
This proves Condition (\ref{l1a}) for this case. Since $\kappa''$ is the only $\vdash_M$-successor of $\kappa$, Lemma~\ref{l1}.\ref{l1b} is satisfied for $\kappa'=\kappa''$ and $\lambda'=\lambda''$.

\medskip
{\boldmath${\rm Stm}_M(l)\in$} $\{(i,{\rm left},l'),\ (i:= a,l'),\ (i, {\rm if}\ a\ {\rm then} \ l',\ {\rm else}\ l'')\}$:\\
In these cases the argument is the same as in the first case.

\medskip

{\boldmath${\rm Stm}_M(l)=(i:=f(i_1,\ldots, i_n), l'):$}\\
Then ${\rm Stm}_N(l)=(i:=g(i_1,\ldots, i_n), l')$ such that
$f : X_{i_1}\times\ldots\times X_{i_n}\mto X_i$ is a $(\gamma_{i_1},\ldots, \gamma_{i_n}, \gamma_i)$-realization of $g : Y_{i_1}\times\ldots\times Y_{i_n}\mto Y_i$.
Since $\lambda$ has a successor, $g(\widetilde y)\neq\emptyset$.
Since $f$ realizes $g$ and $\kappa$ realizes $\lambda$, $\widetilde x$ realizes $\widetilde y$, hence $f(\widetilde x)\neq\emptyset$ by (\ref{f1}).
Since $f(\widetilde x)\neq\emptyset$, $\kappa$ has a successor
by Definition~\ref{d2}.\ref{d2c5}. This proves Condition (\ref{l1a}) for this case lemma.

Let $\kappa'$ be a successor of $\kappa$. Then by Definition~\ref{d2}.\ref{d2c5},
$\kappa'=\kappa[ {\rm cell}_i\leftarrow x_0]\,[{\rm label}\leftarrow l']$ for some
$x_0\in f(\widetilde x)$.
Since $f$ realizes $g$, $\widetilde x$ realizes $\widetilde y$ and $f(\widetilde x)\neq\emptyset$, there is some
$y_0\in \gamma_i(x_0)\cap g(\widetilde y)$ by  (\ref{f1}).
Let $\lambda':=\lambda [ {\rm cell}_i\leftarrow y_0]\,[{\rm label}\leftarrow l']$.
Since $x_0$ realizes $y_0$, $\kappa'$ realizes $\lambda'$. And since
$y_0\in  g(\widetilde y)$, $\lambda\vdash _N\lambda'$ by Definition~~\ref{d2}.\ref{d2c5}.
This proves Condition (\ref{l1b}) for this case.
\medskip

{\boldmath${\rm Stm}_M(l)=({\rm if}\ f(i_1,\ldots,i_n) \ {\rm then}\  l', \ {\rm else}\ l''):$}\\
Then ${\rm Stm}_N(l)=({\rm if}\ g(i_1,\ldots,i_n) \ {\rm then}\  l', \ {\rm else}\ l'')$ such that $f\pf X_{i_1}\times\ldots\times X_{i_n}\to\s$ is a $(\gamma_{i_1},\ldots, \gamma_{i_n}, \Is)$-realization of $g\pf Y_{i_1}\times\ldots\times Y_{i_n}\to\s$.
Since $\lambda$ has a successor, either
$g(\widetilde{y})=0$ or $g(\widetilde{y})=1$. Since $f$ is a $(\gamma_{i_1},\ldots, \gamma_{i_n}, \Is)$-realization of
$g$ and $\kappa$ realizes $\lambda$, $\widetilde{x}$ realizes
$\widetilde{y}$ and either $f(\widetilde{x})=0$ or
$f(\widetilde{x})=1$. By Definition~\ref{d2}.\ref{d2c7} $\kappa$
has a successor. This proves Condition (\ref{l1a}) for this case.

Let $\kappa\vdash_M\kappa'$. Then by Definition~\ref{d2}.\ref{d2c7}
$\kappa'=\kappa[{\rm label}\leftarrow l']$ if $f(\widetilde{x})=0$
and $\kappa'=\kappa[{\rm label}\leftarrow l''] $ if
$f(\widetilde{x})=1$. Since $f$ is a $(\gamma_{i_1},\ldots, \gamma_{i_n}, \Is)$-realization of $g$ and $\widetilde{x}$
realizes $\widetilde{y}$, either $g(\widetilde{y})=0$ or
$g(\widetilde{y})=1$. Let $\lambda'=\lambda[{\rm label}\leftarrow
l']$ if $g(\widetilde{y})=0$ and $\lambda'=\lambda[{\rm
label}\leftarrow l''] $ if $g(\widetilde{y})=1$. Since
$\widetilde{x}$ realizes $\widetilde{y}$, $\kappa'$ realizes
$\lambda'$ and $\lambda\vdash\lambda'$ by
Definition~\ref{d2}.\ref{d2c7}. This proves Condition (\ref{l1b}) for this case.
\qq

\noindent We apply Lemma~\ref{l1} to prove Theorem~\ref{t1}.

\pproof (Theorem~\ref{t1})
First we observe that for configurations $\kappa$ of $\bf M$ and $\lambda$ of~$\bf N$,
\begin{eqnarray}\label{f10}
 (\kappa\mbox{ is accepting}\iff  \lambda\mbox{ is accepting),\ \ \ \ if }
 \kappa \mbox{ realizes } \lambda.
\end{eqnarray}
Let $x=(x_1,\ldots,x_k)\in X_1\times\ldots\times X_k$,
$y=(y_1,\ldots,y_k)\in Y_1\times\ldots\times Y_k$ and $y\in
\gamma_1\times\ldots\times\gamma_k(x)\cap \dom(f_N)$.
 Let $\kappa^0:={\rm IC}_M(x_1,\ldots,x_k)$ and
$\lambda^0:={\rm IC}_N(y_1,\ldots,y_k)$. Then $\kappa^0$ realizes $\lambda^0$.
Since $y\in \dom(f_N)$,
\begin{eqnarray}
\label{f13}
&&\mbox{\hspace{-6ex}there is an accepting computation $(\lambda^0,\ldots,\lambda^n)$ and}\\
\label{f14}
&&\mbox{\hspace{-6ex}every maximal computation with first configuration $\lambda^0$
is accepting.}
\end{eqnarray}

\noindent By Definitions~\ref{d3} and \ref{d2}.\ref{d2e} it suffices to prove:
\begin{enumerate}[(1)]
\item \label{pr3}there is an accepting computation on $\bf M$ with first configuration $\kappa^0$, and
\item \label{pr4}every maximal computation on $\bf M$ with first configuration $\kappa^0$ is an accepting computation $(\kappa^0,\ldots,\kappa^n)$
such that there is an accepting computation $(\lambda^0,\ldots,\lambda^n)$ such that $\kappa^n$ realizes $\lambda^n$.
\end{enumerate}

\noindent Proof of (\ref{pr3}): We know that $\kappa^0$ realizes $\lambda^0$.
For induction suppose  $(\kappa^0,\ldots,\kappa^m)$ is a computation on $\bf M$ and $(\lambda^0,\ldots,\lambda^m)$ is a computation on $\bf N$ such that $\kappa^m$ realizes $\lambda^m$. Suppose $\lambda^m$ has a successor.
By Lemma~\ref{l1}.\ref{l1a} $\kappa^m$ has a successor $\kappa^{m+1}$ and by
Lemma~\ref{l1}.\ref{l1b} there is some successor $\lambda^{m+1}$ of $\lambda^m$ such that $\kappa^{m+1}$ realizes $\lambda^{m+1}$. Then $(\kappa^0,\ldots,\kappa^{m+1})$ and $(\lambda^0,\ldots,\lambda^{m+1})$ are computations such that $\kappa^{m+1}$ realizes $\lambda^{m+1}$. By  (\ref{f14}) this inductive process must end with an accepting computation $(\lambda^0,\ldots,\lambda^n)$.
 For the corresponding computation
$(\kappa^0,\ldots,\kappa^n)$ on $\bf M$, $\kappa^n$ realizes $\lambda^n$. By (\ref{f10}), this computation is accepting.
\smallskip

Proof of (\ref{pr4}):
Let $(\kappa^0,\kappa^1,\ldots$) be a maximal computation on $M$. Then $\kappa^0$ realizes $\lambda^0$.
Assume, for the computation $(\kappa^0,\ldots,\kappa^m)$ we have determined a computation $(\lambda^0,\ldots,\lambda^m)$  on $\bf N$ such that $\kappa^m$ realizes $\lambda^m$.


Suppose, $\lambda^m$ has a successor.  Then, by Lemma~\ref{l1}, $\kappa^m$
has a successor as well. Therefore, $(\kappa^0,\ldots,\kappa^m)$ is not
maximal, hence $\kappa^{m+1}$ exists. By Lemma~\ref{l1}.\ref{l1b},
$\lambda^m$ has a successor $\lambda^{m+1}$ such that $\kappa^{m+1}$
realizes $\lambda^{m+1}$.

By (\ref{f14}) this process must stop at some $n$ such that $(\kappa^0,\ldots,\kappa^n)$ is an initial part of our maximal computation and
$(\lambda^0,\ldots,\lambda^n)$ is accepting. Since $\kappa^n$ realizes $\lambda^n$,
$\kappa^n$ is accepting by (\ref{f10}). This proves~\ref{pr4}.
\qq

\section{Computable Functions on \texorpdfstring{$\s$}{s} and 
  \texorpdfstring{$\om$}{om} are Closed Under Programming}\label{sece}

Suppose, in Definition~\ref{d1}, $X_i=\om$ for all $i$ and all the
functions $f$ in Definitions~\ref{d1}.\ref{d1e5}
and~\ref{d1}.\ref{d1e7} are computable. We want to show that $f_M\pf
(\om)^k\to\om$ is computable. We solve the problem by reduction to
generating functions and sets on $\s$.

Computable functions $f\pf (\om)^k\to \s$ or $f\pf (\om)^k\to \om$
can be generated by monotone computable word functions \cite[Def~2.1.10, Lemma~2.1.11]{Wei00}. Here we use the slightly modified Definition 2 from \cite{Wei08}.

\begin{defi}\label{d4}\hfill
 \begin{enumerate}[(1)]
 \item\label{d4b}
 Call a function $ h\pf (\s)^k\to \s$ {\em monotone-constant}, iff
\\[-2ex]
\[(h(y)\downarrow \ \mbox{and}\ \  y\prf y')
 \imp ( h(y')\downarrow \ \mbox{and}\ \  h(y)=h(y'))\,.\]
For a monotone-constant  function $h$ define $T_*(h)\pf(\om)^k\to\s$ by
\begin{eqnarray}\label{f8}
T_*(h)(x)=w &:\iff& (\exists y\in (\s)^k)\,(y\prf x\an h(y)=w).
\end{eqnarray}

 \item\label{d4c}
Call a function $ h\pf (\s)^k\to \s$ {\em monotone}, iff
\\[-2ex]
\[( h(y)\downarrow  \ \mbox{and}\ \  y\prf y')
\imp (h(y')\downarrow \ \mbox{and}\ \  h(y)\prf h(y'))\,.\]
For a monotone function $h$ define  $T_\omega(h)\pf(\om)^k\to\om$ by
\vspace*{-1ex}
\begin{eqnarray}\label{f9}
T_\omega(h)(x)=q  &:\iff& q=\sup\!{_{_\prf}}\{h(y)\mid  y\prf x
\mbox{ and } h(y)\downarrow\}\,.
\end{eqnarray}
\end{enumerate}
\end{defi}

\noindent Notice that $T_*(h)$ and $T_\omega(h)$ are well-defined by
the ``generating function''~ $h$. By Lemma~\ref{l5}, Turing computable
functions $f\pf (\om)^k\to\s$ or $f'\pf ({\om})^k\to\om$ can be
generated by computable word functions $h\pf (\s)^k\to\s$ which are
monotone-constant or monotone, respectively. We include the continuous
versions.

\begin{lem}\label{l5}\cite{Wei00,Wei08}
\begin{enumerate}[\em(1)]
\item \label{l5a} A function $f\pf (\om)^k\to\s$
is continuous with open domain, iff $f=T_*(h)$ for some monotone-constant
 function $h\pf(\s)^k\to\s$.
\item \label{l5b} A function $f\pf (\om)^k\to\s$
is Turing computable, iff $f=T_*(h)$ for some Turing computable monotone-constant
 function
$h\pf(\s)^k\to\s$.
\item\label{l5c}
A function $f\pf (\om)^k\to\om$ is continuous with $G_\delta$-domain,
iff $f=T_\omega (h)$ for some monotone function $h \pf(\s)^k\to\s$.
\item \label{l5d}
A function $f\pf (\om)^k\to\om$ is Turing computable, iff $f=T_\omega (h)$
for some Turing computable monotone function $h \pf(\s)^k\to\s$.\qed
\end{enumerate}
\end{lem}

\noindent Properties \ref{l5b} and \ref{l5d} are (essentially) \cite[Lemma~2.1.11]{Wei00}.
The proofs show how Type-2 machines can be converted to ``generating'' Turing machines and conversely. For proving the continuous versions we can use machines with an oracle $B\In\s$.
 For the next proofs we extend the prefix relation
$\prf$ on $\s\cup\om$ straightforwardly to $\s\cup\om\cup\Gamma$ and to configurations of machines operating on the sets $\s$ or $\om$. For $u,v\in\s\cup\om\cup\Gamma$ and configurations
$\kappa=(l,(\alpha_0,m_0),\ldots ,(\alpha_L,m_L))$ and $\kappa'=(l',(\alpha'_0,m'_0),\ldots ,(\alpha'_L,m'_L))$
of generalized Turing machines $\bf M$ and $\bf N$, respectively, define:
\[u\prf_1 v :\iff u=v\in\Gamma \ \mbox{or}\  (u,v\in\s\cup\om \ \mbox{and}\  u\prf v)\,,\]
\[\kappa \prf_2 \kappa' :\iff (\forall\;1\leq i\leq L)(\forall\;j\in \IZ)\,)(l=l',\ m_i= m'_i,\
\alpha_i(j)\prf_1  \alpha'_i(j)).\]

\begin{lem}\label{l3}
Let ${\bf M}=({\cal L},l_0,l_f,\Gamma, b, k,L,(X_i)_{0\leq i\leq L},{\rm Stm}_M)$
be a generalized Turing machine such that $X_i=\s$ for $0\leq i\leq L$.
For every $l\in{\cal L}$ let $f\pf(\s)^n\to\s$ be monotone if $\;{\rm Stm}_M(l)=(i:=f(i_1,\ldots, i_n), l')$ and let $f\pf(\s)^n\to\s$ be monotone constant if $\;{\rm Stm}_M(l)=({\rm if}\ f(i_1,\ldots,i_n) \ {\rm then}\  l', \ {\rm else}\ l'')$.
 Then $f_M\pf(\s)^k\to\s$ is monotone.
\end{lem}

\pproof
Since all the functions used in $\bf M$ are single-valued, the successor relation on configurations is a partial function, which we denote by $S$.
First, we prove that for all configurations $\kappa,\kappa'$ of $\bf M$:
\begin{eqnarray}\label{f16}
(S(\kappa)\downarrow\an \kappa\prf_2 \kappa') &\Longrightarrow &
(S(\kappa')\downarrow\an (S(\kappa)\prf_2 S(\kappa'))\,.
\end{eqnarray}
If $\kappa\prf_2 \kappa'$ then $\kappa$ and $\kappa'$ have the same labels and the same head positions. Therefore, they can be written as
$\kappa=(l,(\alpha_0,m_0),\ldots ,(\alpha_L,m_L))$ and
$\kappa'=(l,(\alpha'_0,m_0),\ldots ,(\alpha'_L,m_L))$ such that $\alpha_i(j)=\alpha'_i(j)\in\Gamma$ or $\alpha_i(j)\prf\alpha'_i(j)\in \s$ for all $i,j$. The successor function $S$ changes $\kappa$ and $\kappa'$ only locally. We consider the six alternatives from Definition~\ref{d2}.

\medskip{\boldmath ${\rm Stm}_M(l)=(i,{\rm right},l')$:}\\
Then $S(\kappa)= \kappa[{\rm head}_i\leftarrow m_i+1]\,[{\rm label}\leftarrow l']$
and $S(\kappa')= \kappa'[{\rm head}_i\leftarrow m_i+1]\,[{\rm label}\leftarrow l']$.
Obviously $S(\kappa)\prf_2 S(\kappa')$.

\medskip
{\boldmath ${\rm Stm}_M(l)\in$} $\{(i,{\rm left},l'),\ (i:= a,l'),\ (i, {\rm if}\ a\ {\rm then} \ l',\ {\rm else}\ l'')\}$:\\
In these cases the argument is the same as in the first case.

{\boldmath ${\rm Stm}_M(l)=  (i:=f(i_1,\ldots,i_n),l')  $:}\\
The statement can change only the label and the inscription under the head of Tape $i$.
Since $S(\kappa)$ exists, $x:=f(x_{i_1},\ldots,x_{i_n})$ exists, where $x_{i_j}=\alpha_j(m_j)\in\s$ for $1\leq j\leq n$ (see Definition~\ref{d2}). Since $\kappa\prf_2\kappa'$,
$x_{i_j}\prf x'_{i_j}:=\alpha'_j(m_j)\in\s$ ($1\leq j\leq n$). Since $f$ is monotone,
$f(x'_{i_1},\ldots,x'_{i_n})$ exists and $x=f(x_{i_1},\ldots,x_{i_n})\prf f(x'_{i_1},\ldots,x'_{i_n})= x'$.
Since $\kappa\prf_2\kappa'$ and
$S(\kappa)=\kappa[ {\rm cell}_i\leftarrow x]\,[{\rm label}\leftarrow l']$
and $S(\kappa')=\kappa'[ {\rm cell}_i\leftarrow x']\,[{\rm label}\leftarrow l']$,
$S(\kappa)\prf_2 S(\kappa')$.

\medskip
{\boldmath ${\rm Stm}_M(l)=  ({\rm if}\ f(i_1,\ldots,i_n)\ {\rm then}\ l',{\rm else}\ l'')  $:}\\
The statement changes only the labels (or cannot be applied). Since $S(\kappa)$ exists by assumption, $f(x_{i_1},\ldots,x_{i_n})\in\{0,1\}\In\s$ exists, where $x_{i_j}=\alpha_j(m_j)\in\s$ for $1\leq j\leq n$ (see Definition~\ref{d2}).
Since $\kappa\prf_2\kappa'$,
$x_{i_j}\prf x'_{i_j}:=\alpha'_j(m_j)\in\s$ ($1\leq j\leq n$).
Since $f$ is monotone constant,
$f(x'_{i_1},\ldots,x'_{i_n})$ exists and $f(x_{i_1},\ldots,x_{i_n})= f(x'_{i_1},\ldots,x'_{i_n})$.
By Definition~\ref{d2}, also the labels of $S(\kappa)$ and $S(\kappa')$ are the same,
hence $S(\kappa)\prf_2 S(\kappa')$.

\medskip
\noindent This proves (\ref{f16}).

Now suppose $w=(w_1,\ldots,w_k)\prf (w'_1,\ldots,w'_k)=w'$.
Then ${\rm IC}_M(w)\prf_2{\rm IC}_M(w')$ (Definition~\ref{d2}).
Suppose $f_M(w)$ exists. Then for some $n$, $S^n\circ {\rm IC}_M(w)$ exists and is an accepting configuration such that $f_M(w)=\alpha_0(0)$. From (\ref{f16}) by induction
for all $k\leq n$, $S^k\circ {\rm IC}_M(w')$ exists and
$S^k\circ {\rm IC}_M(w)\prf_2 S^k\circ {\rm IC}_M(w')$.
Since $S^n\circ {\rm IC}_M(w)\prf_2 S^n\circ {\rm IC}_M(w')$ and $S^n\circ {\rm IC}_M(w)$ is accepting, $S^n\circ {\rm IC}_M(w')$ is accepting.
Then, $f_M(w)$, the word on (Tape $0$, Cell $0$) of
$S^n\circ {\rm IC}_M(w)$ is a prefix of $f_M(w')$, the word on (Tape $0$, Cell $0$) of $S^n\circ {\rm IC}_M(w')$, hence $f_M(w)\prf f_M(w')$.
Therefore, $f_M$ is monotone.
\qq

Let $\bf M$ be a generalized Turing machine on generating word functions and let $\bf N$ be the corresponding generalized Turing machine on generated functions on $\om$. Then
$f_M\pf(\s)^k\to\s$ generates an extension of $f_N\pf(\om)^k\to\om$:

\begin{thm}\label{t2}
Let ${\bf M}=({\cal L},l_0,l_f,\Gamma, b, k,L,(X_i)_{0\leq i\leq L},{\rm Stm}_M)$ and\\
${\bf N}=({\cal L},l_0,l_f,\Gamma, b, k,L,(Y_i)_{0\leq i\leq L},{\rm Stm}_N)$
be generalized Turing machines such that $X_i=\s$ and $Y_i=\om$ for $0\leq i\leq L$
and all functions  occurring in $\bf M$ or $\bf N$ are single-valued.
Assume that for all labels $l\in {\cal L}$,

\begin{enumerate}[\em(1)]
\item \label{t2a}
if $\;{\rm Stm}_M(l)\in\{ (i,{\rm right},l'),\ (i,{\rm left},l'),\ (i:= a,l'),\ (i, {\rm if}\ a\ {\rm then} \ l',\ {\rm else}\ l'')\}$ then
$\;{\rm Stm}_M(l) ={\rm Stm}_N(l)$,

\item \label{t2b} if $\;{\rm Stm}_M(l)=(i:=f(i_1,\ldots, i_n), l')$ then
$f$ is monotone and $\;{\rm Stm}_N(l)=(i:=g(i_1,\ldots, i_n), l')$ such that
$T_\omega(f)$ extends $g$,

\item \label{t2d}
if
$\;{\rm Stm}_M(l)=({\rm if}\ f(i_1,\ldots,i_n) \ {\rm then}\  l', \ {\rm else}\ l'')$   then $f$ is monotone-constant and
 $\;{\rm Stm}_N(l)=({\rm if}\ g(i_1,\ldots,i_n) \ {\rm then}\  l', \ {\rm else}\ l'')$      such that $T_*(f)$ extends $g$.
\end{enumerate}
Then $T_\omega(f_M)$ extends $f_N\pf(\om)^k\to\om$.
\end{thm}

\pproof We must prove that for all  $q\in\dom(f_N)$,
\begin{eqnarray}\label{f11}
f_N(q)=T_\omega(f_M)(q)={\sup_{}}{\phantom{.}}_{\hspace{-.5ex}\prf} \{f_M(u)\mid u\in(\s)^k, u\prf q\ {\rm and}\ f_M(u)\downarrow\}\,.
\end{eqnarray}
Since all the functions used in $\bf M$ and $\bf N$ are single-valued, the successor relations on configurations  are functions, which we denote by $S$ for both machines.
For a word $w\in\s$ let $|w|$ denote its length. For a configuration $\kappa$ for $\bf M$ define the precision by
\[P(\kappa):= \min\{|\alpha_i(j)|\mid 0\leq i\leq L,\ j\in\IZ,\  \alpha_i(j)\in\s\}\,.  \]
For $q=(q_1,\ldots,q_k)\in(\om)^k$ and $e\in\IN$ let
$q^{<e}:=(w_1,\ldots,w_k)$ where $w_i$ is the prefix of $q_i$ of length $e$.

\begin{prop}\label{pr1}Suppose, $q=(q_1,\ldots,q_k)\in\dom(f_N)\In(\om)^k$.
Then for all $m$ such that $\lambda:=S^m\circ {\rm IC}_N(q)$ exists:
\begin{eqnarray}\label{f17}
(\forall d)(\exists \overline e)(\forall e\geq \overline e)
(\kappa:=S^m\circ {\rm IC}_M( q^{<e}  )\downarrow,\  \kappa\prf_2 \lambda, \  P(\kappa) \geq d)
\end{eqnarray}
(where $d,\,\overline e,\,e\in\IN$). This means that for sufficiently precise input,  $\kappa$ exists and approximates $\lambda$ with at least precision~$d$.
\end{prop}

\pproof (Proposition~\ref{pr1})
 We prove (\ref{f17}) by induction on $m\in\IN$.

\medskip
{\boldmath $m=0$:} For $d\in\IN$ choose $\overline e:=d$.
 Then for $e\geq \overline e$ by Definition~\ref{d2}.\ref{d2e},  $S^0\circ {\rm IC}_M(q^{<e})\prf_2
 S^0\circ {\rm IC}_N(q)$ and $P( S^0\circ{\rm IC}_M(q^{<e}))=e\geq d$.

\medskip
{\boldmath $m \Longrightarrow m+1$:} Assume that the statement has been proved for $m$ and assume that $S^{m+1}\circ {\rm IC}_N(q)$ exists. Then $S^m\circ {\rm IC}_N(q)$ exists and can be written as
\begin{eqnarray}\label{f12}
\lambda &:=& S^m\circ {\rm IC}_N(q)=(l,(\beta_0,m_0),\ldots,(\beta_L,m_L))\,.
\end{eqnarray}
Let $d\in\IN$. We consider the 6 alternatives  for ${\rm Stm}_N(l)$ from Definition~\ref{d2}. Notice that the successor functions $S$ of $\bf M$ and $\bf N$ change configurations only locally.

{\boldmath ${\rm Stm}_N(l)=(i,{\rm right},l')$:}\\
By assumption, there is some $\overline e\in\IN$ such that for all $e\geq \overline e$,
$\kappa:=S^m\circ {\rm IC}_M( q^{<e}  )$ exists, $\kappa\prf_2 \lambda$ and $  P(\kappa) \geq d $.
We show that we can choose this number $\overline e$ for $m+1$ and $d$ as well.
Since $\kappa\prf_2  \lambda$, $\kappa$ can be written as
\[\kappa=S^m\circ {\rm IC}_M( q^{<e} )=(l,(\alpha_0,m_0),\ldots,(\alpha_L,m_L))\]
such that for all $i$ and $j$, $\alpha_i(j)=\beta_i(j)\in\Gamma$ or
  $\alpha_i(j)\prf\beta_i(j)$ (where $\alpha_i(j)\in \s$ and  $\beta_i(j)\in\om$).
By the condition in Theorem~\ref{t2}.\ref{t2a},  ${\rm Stm}_M(l)={\rm Stm}_N(l)=(i,{\rm right},l')$. Therefore,
\[\begin{array}{lll}S(\lambda )&=&\lambda[{\rm head}_i\leftarrow m_i+1]\,[{\rm label}\leftarrow l']\,, \\
S(\kappa )&=&\kappa[{\rm head}_i\leftarrow m_i+1]\,[{\rm label}\leftarrow l']\,.
\end{array}\]
Since on $\lambda$ and $\kappa$ the successors $S$ operate in the same way depending at most on tape cells containing elements of $\Gamma$ and changing at most such tape cells,  $\kappa\prf_2\lambda$ implies $S\circ \kappa\prf_2 S\circ \lambda$
and $P(S\circ \kappa)=P(\kappa)\geq d$.

\medskip
{\boldmath ${\rm Stm}_N(l)\in\{(i,{\rm left},l'),\    (i:=a,l') \  (i,{\rm if}\ a\ {\rm then}\ l',{\rm else}\ l'')\}$:}\\
The arguments in these cases are the same as in the first case.\\

{\boldmath ${\rm Stm}_N(l)=(i:=g(i_1,\ldots,i_n),l')$:}\\
By the condition in Theorem~\ref{t2}.\ref{t2b},  ${\rm Stm}_M(l)=(i:=f(i_1,\ldots,i_n),l')$ such that $T_\omega(f)$ extends $g$.

Let $s:=(s_1,\ldots,s_n)\in(\om)^n$ such that $s_j=\beta_{i_j}(m_{i_j})$ is the content of the cell under the head of  Tape $i_j$ of the configuration $\lambda$ (see (\ref{f12})). Since
$S^{m+1}\circ {\rm IC}_N(q)=S(\lambda)$ exists, $s\in\dom(g)$.
 Since $T_\omega(f)$ extends $g$, by the sup-condition in Definition~\ref{d4}.\ref{d4c}
\begin{eqnarray}\label{f18}
(\exists \overline b)(\forall b\geq \overline b)
(f(s^{<b})\downarrow,\ \ f(s^{<b})\prf g(s)\ \mbox{and}\ \ |f(s^{<b})|\geq d)\,.
\end{eqnarray}
By induction as  a special case of (\ref{f17}), for $\max(\overline b,d)$ there is some $\overline e$ such that
\begin{eqnarray}\label{f19}
(\forall e\geq \overline e)
(\kappa:=S^m\circ {\rm IC}_M( q^{<e})\downarrow,\ \ \kappa\prf_2 \lambda
\ \ \mbox{and}\ \ P(\kappa) \geq \max (\overline b,d))\,.
\end{eqnarray}

We show that this constant $\overline e$ is appropriate for $m+1$ and $d$ in (\ref{f17}).
Let $e\geq \overline e$, $\kappa:=S^m\circ {\rm IC}_M( q^{<e}  )$ and $\kappa':= S^{m+1}\circ {\rm IC}_M( q^{<e})=S(\kappa)$.
Since $\kappa\prf_2  \lambda$, $\kappa$ can be written as (see~(\ref{f12}))
\[\kappa=S^m\circ {\rm IC}_M( q^{<e} )=(l,(\alpha_0,m_0),\ldots,(\alpha_L,m_L))\]
such that for all $0\leq j\leq L$ and $i'$, $\alpha_j(i')=\beta_j(i')\in\Gamma$ or
  $\alpha_j(i')\prf\beta_j(i')$ (where $\alpha_j(i')\in \s$ and  $\beta_j(i')\in\om$). Let $u:=(u_1,\ldots,u_k)$ where $u_j:=\alpha_{i_j}(m_{i_j})$. Then $u\prf s$.
By Definition~\ref{d2}
\begin{eqnarray}\label{f20}
\begin{array}{lll}
S(\kappa )&=&\kappa[{\rm cell}_i\leftarrow v]\,[{\rm label}\leftarrow l']\,,\\
S(\lambda )&=&\lambda[{\rm cell}_i\leftarrow q']\,[{\rm label}\leftarrow l']\,,
\end{array}
\end{eqnarray}
where $v=f(u)$ and $q'=g(s)$.  By (\ref{f17}) we must prove:
\begin{eqnarray}\label{f22}
\kappa':=S^{m+1}\circ {\rm IC}_M( q^{<e}  )\downarrow,\  \kappa'\prf_2 S(\lambda)\ \mbox{and} \  P(\kappa') \geq d\,.
\end{eqnarray}
 Since $P(\kappa)\geq \overline b$,
$s^{<\overline b}\prf u$. Since $f(s^{<\overline b})\downarrow$ by (\ref{f18}) and $f$ is monotone, $f(u)$ exists. Therefore, $\kappa'=S(\kappa)$ exists.
Then $S(\kappa)$ and $S(\lambda)$ can be written as
\begin{eqnarray}\label{f21}
\begin{array}{lll}
S(\kappa) & = & (l',(\alpha'_0,m_0),\ldots,(\alpha'_L,m_L))\,,\\
S(\lambda) & = &(l',(\beta'_0,m_0),\ldots,(\beta'_L,m_L))\,.
\end{array}
\end{eqnarray}
$S(\kappa)$ differs from $\kappa$ only on Cell $i$, the cell under the head of Tape $i$, and $S(\lambda)$ differs from $\lambda$ only on Cell $i$, the cell under the head of Tape $i$.

Since $u\prf s$ and $T_\omega(f)$ extends $g$, by (\ref{f20}), $\alpha'_i(m_i)=f(u)\prf g(s)= \beta'_i(m_i)\,.$
For all  $(j,i')\neq (i,m_i)$ by $\kappa\prf_2\lambda$ (\ref{f19})
$\alpha'_j(i')= \alpha_j(i')\prf_1 \beta_j(i')= \beta'_j(i')\,.$
Therefore, $S(\kappa)\prf_2 S(\lambda)$.

By (\ref{f19}), $|u_j|\geq \overline b$ for $1\leq j\leq n$, hence $s^{<\overline b}\prf u$, since $u\prf s$. By (\ref{f18}) and monotonicity of $f$, $|\alpha'_i(m_i)|=|f(u)|\geq |f(s^{<\overline b})|\geq d$.
For all  $(j,i')\neq (i,m_i)$ such that $\alpha'_j(i')\in\s$,
$|\alpha'_j(i')|= |\alpha_j(i')| \geq d$, since $P(\kappa)\geq d$ by (\ref{f19} ). Therefore, $P(\kappa')\geq d$.

This proves (\ref{f22}) and finishes the case
${\rm Stm}_N(l)=(i:=g(i_1,\ldots,i_n),l')$.
\medskip

{\boldmath ${\rm Stm}_N(l)=({\rm if}\ g(i_1,\ldots,i_n)\ {\rm then}\ l',{\rm else}\ l'')$:}\\
The proof can be obtained by straightforward modification of the proof of the previous case.
\hfill $\Box$(Proposition~\ref{pr1})
\medskip

It remains to prove (\ref{f11}).
Suppose, $f_N(q)$ exists. Then for some $m\in \IN$,  $\lambda :=S^m\circ {\rm IC}_N(q)$ exists, $\lambda $ is a final configuration and ${\rm OC}(\lambda)=f_N(q)$. Let $d\in\IN$. By Proposition~\ref{pr1} there is some $\overline e\in\IN$ such that
$\kappa:=S^m\circ {\rm IC}_M( q^{<\overline e})$ exists, $\kappa\prf_2 \lambda$ and   $P(\kappa) \geq d$.
Since $\kappa\prf_2\lambda$, also $\kappa$ is a final configuration and $f_M( q^{<\overline e})={\rm OC}(\kappa)\prf {\rm OC}(\lambda)=f_N(q)$.
Furthermore, $|f_M( q^{<\overline e})|\geq d$, since $P(\kappa) \geq d$.
Since $f_M$ is monotone by Lemma~\ref{l3}, $f_N(q)=\sup _\prf\{f_M(q^{<e})\mid e\in\IN\}$.
Since for all $u\in (\s)^k$ with $u\prf q$ there is some $e$ such that $u\prf q^{<e}$,
$\sup _\prf\{f_M(q^{<e})\mid e\in\IN\}=\sup _\prf\{f_M(u)\mid u\prf q\}=T_\omega(f_M)$.

Therefore, $T_\omega(f_M)$ extends $f_N\pf(\om)^k\to\om$.
\qq

If all functions on $\s$ used in the machine $\bf M$ from Theorem~\ref{t2} are computable, then $f_M$ is a computable word function.

\begin{lem}\label{l2}
Let ${\bf M}=({\cal L},l_0,l_f,\Gamma, b, k,L,(X_i)_{0\leq i\leq L},{\rm Stm}_M)$
be a generalized Turing machine such that $X_i=\s$ for $0\leq i\leq L$ and all functions on $\s$ used in the machine are computable. Then $f_M\pf(\s)^k\to\s$ is computable.
\end{lem}

\pproof From the generalized Turing machine $\bf M$ an ordinary Turing machine $\bf N$ computing $f_M$ can be constructed by standard techniques.
\qq

By the next theorem the continuous as well as the computable functions on $\om$ are closed under programming.

\begin{thm}\label{t3}
Let ${\bf N}=({\cal L},l_0,l_f,\Gamma, b, k,L,(Y_i)_{0\leq i\leq L},{\rm Stm}_N)$
be a generalized Turing machine on $\om$, that is, $Y_i=\om$ for $0\leq i\leq L$. If all the functions on $\om$ used in the machine $\bf N$
\begin{enumerate}[\em(1)]
\item \label{t3a}are continuous, then $f_N\pf(\om)^k\to\om$ is continuous,
\item \label{t3b}are computable, then $f_N\pf(\om)^k\to\om$ is computable.
\end{enumerate}
\end{thm}

\pproof (\ref{t3a}) Every function used in $\bf N$ is generated by a monotone or monotone constant word function (Lemma~\ref{l5}). Let $\bf M$ be the machine constructed with these word functions that satisfies the conditions from Theorem~\ref{t2}. Then $T_\omega(f_M)$ extends $f_N$. By Lemma~\ref{l5}, $f_N$ is continuous.

(\ref{t3b}) In addition to Case~(\ref{t3a}) there are even computable word functions. Again the function $f_M$ generates $f_N$. By Lemma~\ref{l2}, $f_M$ is computable, hence $f_N$ is computable by Lemma~\ref{l5}.
\qq

The generalization from $\om$ to $\s$ and $\om$ is straightforward.

\begin{cor}\label{co2}
Theorem~\ref{t3} holds accordingly, if $Y_i\in\{\s,\om\}$ for $0\leq i\leq L$.
\end{cor}

\pproof Define a standard representation $\beta\pf \om\to\s$ of $\s$ by $\beta(\iota(w)0^\omega):=w$ (where $\iota(a_1\ldots a_n):=110a_10\ldots0a_n011$,\cite[Definition~2.1.7]{Wei00}). For $Y_j$ let $\delta_j:=\id_\om$ if $Y_j=\om$ and $\delta_j:=\beta$ if $Y_j=\s$. Then a function $f\pf Y_{i_1}\times \ldots\times Y_{i_n}\to Y_{i_0}$ is computable, iff it is $(\delta_{i_1},\ldots,\delta_{i_n},\delta_{i_0})$-computable by a realization on $\om$.
Let $\bf M$ be a machine obtained from $\bf N$ by replacing every function $f$ on $\om$ and $\s$ by a computable realizing function on $\om$. Then $\bf M$ realizes $\bf N$, hence $f_M$ realizes $f_N$ by Theorem~\ref{t1}. Since $f_M$ is computable by Theorem~\ref{t3}, $f_N$ is computable.
For ``continuous'' the argument is the same.
\qq

\section{Machines on Represented Sets, the Main Results}\label{secf}

After the preparations in Sections~\ref{secd} and \ref{sece} we can easily prove our main results, Theorems~\ref{t4} and~\ref{t5}.
Since the computable functions on $\om$ are closed under composition, the composition of computable functions on represented sets is computable \cite[Theorem~3.1.6]{Wei00}.
The following main result of this article generalizes this observation from single-valued to multi-valued functions and representations and from composition to generalized Turing machines.

\begin{defi}\label{d6}
Let ${\bf P}=({\cal L},l_0,l_f,\Gamma, b, k,L,(Z_i)_{0\leq i\leq L},{\rm Stm}_P)$
be a generalized Turing machine and for each $i$, $0\leq i\leq L$, let $\delta_i:\om\mto Z_i$ be a multi-representation.
 The machine is called $(\delta_i)_{0\leq i\leq L}$-computable, iff
\begin{enumerate}[(1)]
\item \label{d6a} for every statement ``$(i:=f(i_1,\ldots,i_n),l')$'' in $\bf P$ \\
the multi-function $f$
 is $(\delta_{i_1},\ldots,\delta_{i_n},\delta_i)$-computable  and
\item \label{d6b} for every statement  ``$(\mbox{\rm if } f(i_1,\ldots,i_n) \mbox{ \rm then } l', \mbox{ \rm else } l'')$'' in $\bf P$ \\
the  partial function $f$  is $(\delta_{i_1},\ldots,\delta_{i_n},\Is)$-computable.
\end{enumerate}
 ``$(\delta_i)_{0\leq i\leq L}$-continuous'' is defined in the same way with ``continuous'' replacing ``computable''.
\end{defi}

\begin{thm}\label{t4}
Let ${\bf P}$
be a generalized Turing machine with $k$ input tapes and for $0\leq i\leq L$ let $\delta_i$ be a multi-representation of $Z_i$ such that the machine is $(\delta_i)_{0\leq i\leq L}$-computable. Then the function $f_P$ computed by the machine is $(\delta_1,\ldots,\delta_k,\delta_0)$-computable.

Correspondingly with ``continuous'' instead of ``computable''.
\end{thm}

\pproof
Case ``continuous'':
 There is a generalized Turing machine $\bf M$ on $\om$ containing only continuous functions that realizes $\bf  P$ via  $(\delta_i)_{0\leq i\leq L}$ (replace every function in $\bf P$ by a realizing function on $\om$). By Theorem~\ref{t1}, $f_M$ realizes $f_P$ via $(\delta_1,\ldots,\delta_k,\delta_0)$.
By Theorem~\ref{t3}, $f_M$ is continuous. Therefore, $f_P$ is $(\delta_1,\ldots,\delta_k,\delta_0)$-continuous.

The case ``computable'' can be  proved in the same way.
\qq

\begin{cor}\label{co1}
Theorem~\ref{t4} remains true if in Definition~\ref{d6} and in the theorem some
multi-representations $\delta_j:\om\mto Z_j$ are replaced by multi-notations $\nu_j:\s\mto Z_j$.
\end{cor}

\pproof For every multi-notation $\nu:\s\mto X$ there is a multi-representation $\delta:\om\mto X$ such that $\nu\equiv \delta$, and equivalent multi-notations/representations of a set $X$ induce the same computability and continuity on $X$ \cite[Section~8]{Wei08}. Replace every multi-notation by an equivalent multi-representation and apply Theorem~\ref{t4}. In the final result return to the multi-notations.
\qq

In applications, generalized representations are already used informally, whenever defining explicitly Type-2 Turing machines for realizing functions on $\om$ is too cumbersome.

\begin{exa}\label{ex1} \rm  Let $\rm SRI$ be the set of all sequences of intervals $[a;b]\In \IR$ with rational end points $a <b$.
Let $\gamma\pf \om\to {\rm SRI}$ be a canonical representation and define a generalized representation $\delta\pf {\rm SRI}\to \IR$ of the real numbers by
$\delta(I_0,I_1,\ldots)=x\iff \{x\}=\bigcap _n I_n$. Then $\delta\circ\gamma\pf \om\to\IR$ is a representation that is equivalent to the standard representation $\rho$ of the real numbers \cite[Chapter~4]{Wei00}. For proving that addition on the real numbers is computable via $\delta\circ \gamma$ consider the following binary function $f_+$ on $\rm SRI$:
\[ f_+((I_0,I_1,\ldots),(J_0,J_1,\ldots)):=(I_0+J_0,I_1+J_1,\ldots).\]
The experienced reader knows that the function $f_+$ is $(\gamma,\gamma,\gamma)$-computable and a simple proof shows that $f_+$ is a
$(\delta,\delta,\delta)$-realization of addition. By the next lemma from \cite{Wei08} we may conclude that addition is computable via $\delta\circ \gamma$.
\qq
\end{exa}

For multi-functions $\gamma:X\mto Y$ and $\delta: Y\mto Z$  the ``relational'' composition $\delta\odot \gamma$ is defined by
\begin{eqnarray}\label{f23}
z\in\delta\odot\gamma(x)\iff (\exists\, y)\,(y\in\gamma(x)\an z\in\delta(y)),
\end{eqnarray}
see \cite[Sections~3 and 6]{Wei08}. If $\gamma $ and $\delta$ are multi-representations, an element $x\in X$ should be considered as a name of $z$ via the combination of $\gamma$ and $\delta$, if there is some $y\in Y$ such that $x$ is a $\gamma$-name of $y$ and $y$ is a $\delta$-name of $z$, that is, $y\in\gamma(x)$ and $z \in\delta(y)$, hence $z\in(\delta\odot\gamma)(x)$.
Therefore, we use relational composition for multi-representations. Notice that for single-valued $\gamma$ (as in Example~\ref{ex1}), $\delta\odot\gamma=\delta\circ\gamma$.

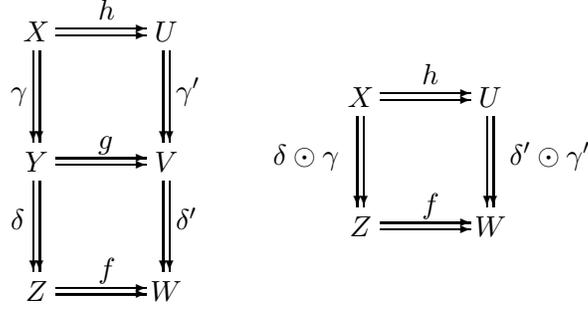
\begin{figure}[ht]

\unitlength0.26ex
\begin{picture}(40,75)(70,0)

\put(0,80){\makebox(2,2){$X$}}
\put(0,40){\makebox(2,2){$Y$}}
\put(0,0) {\makebox(2,2){$Z$}}
\put(40,80){\makebox(2,2){$U$}}
\put(40,40){\makebox(2,2){$V$}}
\put(40,0){\makebox(2,2){$W$}}

\put(7,0){\vector(1,0){28}}
\put(7,40){\vector(1,0){28}}
\put(7,80){\vector(1,0){28}}

\put(7,2){\vector(1,0){28}}
\put(7,42){\vector(1,0){28}}
\put(7,82){\vector(1,0){28}}

\put(0,35){\vector(0,-1){28}}
\put(40,35){\vector(0,-1){28}}
\put(0,75){\vector(0,-1){28}}
\put(40,75){\vector(0,-1){28}}

\put(2,35){\vector(0,-1){28}}
\put(42,35){\vector(0,-1){28}}
\put(2,75){\vector(0,-1){28}}
\put(42,75){\vector(0,-1){28}}

\put(20,85){$h$}
\put(20,45){$g$}
\put(20,5){$f$}

\put(-7,20){$\delta$}
\put(-7,60){$\gamma$}
\put(44,20){$\delta'$}
\put(44,60){$\gamma'$}

\put(100,60){\makebox(2,2){$X$}}
\put(100,20){\makebox(2,2){$Z$}}
\put(140,60){\makebox(2,2){$U$}}
\put(140,20){\makebox(2,2){$W$}}

\put(107,20){\vector(1,0){28}}
\put(107,60){\vector(1,0){28}}
\put(107,22){\vector(1,0){28}}
\put(107,62){\vector(1,0){28}}

\put(100,55){\vector(0,-1){28}}
\put(140,55){\vector(0,-1){28}}
\put(102,55){\vector(0,-1){28}}
\put(142,55){\vector(0,-1){28}}

\put(120,65){$h$}
\put(120,25){$f$}

\put(74,40){$\delta\odot\gamma$}
\put(147,40){$\delta'\odot\gamma'$}
\end{picture}

\caption {Realization is downwards transitive.}
\label{bildh}
\end{figure}

Realization is downwards transitive. If $h$ realizes $g$ and $g$ realizes $f$ then $h$
realizes $f$ w.r.t. the composed representations (Figure~\ref{bildh}).

\begin{lem}[\cite{Wei08}]\label{l19}
Let $\gamma : X\mto Y$, $\delta : Y\mto Z$, $\gamma' : U\mto V$ and
$\delta' : V\mto W$ be generalized multi-representations.
If $h : X\mto U$ is a $(\gamma,\gamma')$-realization of $g : Y\mto V$ and
$g : Y\mto V$ is a $(\delta,\delta')$-realization of $f : Z\mto W$, then
$h : X\mto U$ is a $(\delta\odot \gamma,\delta'\odot \gamma')$-realization of $f : Z\mto W$.\qed
\end{lem}

\noindent In Example~\ref{ex1}, $\IR$ can be called the set of
``abstract'' data, $\om$ the set of ``concrete'' data and $\rm SRI$
the set of data of ``intermediate abstraction''. If computability on
data of intermediate abstraction is well understood, it may be of
advantage to use them as names in generalized representations. The
following theorem generalizes Theorem~\ref{t4}. It shows that the
function $f_P$ computed by a generalized Turing machine $\bf P$
containing only functions realized by computable functions on data of
intermediate abstraction is computable.

\begin{thm}\label{t5}
Let ${\bf P}=({\cal L},l_0,l_f,\Gamma, b, k,L,(Z_i)_{0\leq i\leq L},{\rm Stm}_P)$
be a generalized Turing machine. For each $i$, $0\leq i\leq L$, let $\delta_i:Y_i\mto Z_i$ be a generalized multi-representation and let $\gamma_i:\om\mto Y_i$ be a multi-representation. Suppose,
\begin{enumerate}[\em(1)]
\item \label{t5d} for every statement ``$(i:=f(i_1,\ldots,i_n),l')$'' in $\bf P$
the multi-function $f$ has a realization $g$ via
 $(\delta_{i_1},\ldots,\delta_{i_n},\delta_i)$ that is $(\gamma_{i_1},\ldots,\gamma_{i_n},\gamma_i)$-computable, and
\item \label{t5e} for every statement  ``$(\mbox{\rm if } f(i_1,\ldots,i_n) \mbox{ \rm then } l', \mbox{ \rm else } l'')$'' in $\bf P$
the  partial function $f$ has a realization $g$ via $(\delta_{i_1},\ldots,\delta_{i_n},\Is)$ that is $(\gamma_{i_1},\ldots,\gamma_{i_n},\Is)$-computable.
\end{enumerate}
Let $\bf M$ be a machine obtained from $\bf P$ by replacing every function $f$ by some function $g$ computable w.r.t the $\gamma_j$ realizing $f$ via the $\delta_j$ as described in (\ref{t5d}) and (\ref{t5e}). Then
\begin{enumerate}[\em(a)]
\item \label{t5a} $f_M$ is $(\gamma_1,\ldots,\gamma_k,\gamma_0)$-computable,
\item \label{t5b} $f_M$ realizes $f_P$ via $(\delta_1,\ldots,\delta_k,\delta_0)$,
\item \label{t5c} $f_P$ is $(\delta_1\odot\gamma_1,\ldots,\delta_k\odot\gamma_k,\delta_0\odot\gamma_0)$-computable.
\end{enumerate}
\end{thm}

\pproof \hfill
\begin{enumerate}[(a)]
\item This follows from Theorem~\ref{t4}.
\item This follows from Theorem~\ref{t1}.
\item This follows from (a) and (b) by Lemma~\ref{l19}.\qed
\end{enumerate}
\newpage

\section{Examples}\label{secg}
The feasible real RAM \cite{BH98}, a machine model for real computation, allows approximate multi-valued branching.\\
 $\leq_k:\IR\times\IR\mto \{{\rm tt,ff}\}$,\ \ \
$ \leq_k(x,y)\left\{
\begin{array}{ll}
= {\rm tt} & \mbox{if}\ \ x<y\\
\in \{{\rm tt,ff}\}& \mbox{if}\ \ y\leq x\leq y+2^{-k}\\
= {\rm tt} & \mbox{if}\ \ y+2^{-k}<x\,.
\end{array}\right. $\\
This is not allowed in generalized Turing machines but can be simulated by a multi-valued function followed by a single-valued test.

\smallskip
Theorems~\ref{t4} and~\ref{t5} allow to formulate algorithms and argue about them in a more abstract way which is closer to ordinary analysis and which usually is simpler and more transparent. In Example~\ref{ex1} we have used sequences of rational intervals as names of real numbers.

\smallskip
As another example consider $C^\infty(\IR)$ the set of all infinitely often differentiable real functions. There is a canonical representation $\beta$ of $C(\IR)$, then $\gamma:=[\beta]^\omega$ is a canonical representation of $(C(\IR))^\omega$
\cite[Definition~3.3.3]{Wei00}. Define a generalized representation
$\delta\pf (C(\IR))^\omega \to C^\infty(\IR)$ as follows:
\[\delta(f_0,f_1,\ldots)=g\iff (\forall i)\, f_i=g^{(i)}\,\]
(a name of $g$ is a list of all of its derivatives). This generalized representation may be useful in the study of distributions \cite{ZW03}.

\medskip
As another example we consider computing the sum
$s:=\sum_{j=0}^\infty a_jz^j$ of a complex power series. Let $R$ be the radius of convergence and $s_n:=\sum_{j=0}^{n-1} a_jz^j$ the partial sum of the first $n$ terms. Let  $r<R$, $r\in\IQ$, and let $M\in\IQ$ be a Cauchy constant for $r$, that is, $(\forall j)\,|a_j|\leq M\cdot r^{-j}$. Then for all $|z|< r$,
\[|s_n- s|\leq M\frac{(|z|/r)^n}{1-(|z|/r)}\]
We want to show that the operator $H: ((a_j)_j,r,M, z)\to s$ is computable via the standard representations of the occurring sets \cite{Wei00}. In the following we say ``computable'' instead of ``computable via the standard representations''.

First, from the inputs $(a_j)_j$, $r$, $M$, $z$ and $k$ we compute some complex number $b_k\in\IC$ such that $|b_k-s|\leq 2^{-k}$ as follows:\\
-- compute $c:=|z|/r$\\
-- find some $q\in\IQ$ such that $c<q<1$, \ \ (then $|s_n- s|\leq M\cdot q^n/(1-q)$)\\
-- find some $n\in\IN$ such that $M\cdot q^n/(1-q)<2^{-k}$, \ (then $|s_n- s|\leq 2^{-k}$)\\
-- For $m=0,1,\ldots,n$ compute in turn $z^m$ and $s_m$.\\
-- Let $b_k:=s_n$ be the result.\\
It is easy to find a generalized Turing machine $\bf N$ for this algorithm that uses the arithmetical operations on $\IN, \IQ ,\IR$ and $\IC$  and the projection $((a_j)_j,m)\mapsto a_m$ all of which are computable \cite{Wei00}. Therefore, by Theorem~\ref{t4}, $f_N:((a_j)_j,r,M, z,k)\mmto b_k$ is computable. By \cite[Theorem~35]{Wei00} the multi-function $S\circ f_N:((a_j)_j,r,M, z)\mmto (b_k)_{k\in\IN}$ is computable, where
\[(b_k)_k\in S\circ f_N((a_j)_j,r,M, z)\iff
(\forall k)\,b_k\in  f_N((a_j)_j,r,M, z,k)\,.\]
Since $(b_k)_k$ is a sequence of complex numbers such that $|s-b_k|\leq 2^{-k}$ and the limit operator ${\rm Lim}:(b_k)_k\mapsto \lim_{k\to\infty}b_k$ is computable (cf. \cite[Theorem~4.3.7]{Wei00}, $H={\rm Lim}\circ S\circ f_N$ is computable.

Although a multi-function has been used in the determination of $q$, the function $H$ is single-valued. Notice that the only informal part in the above proof is the specification of the generalized Turing machine $\bf N$. But this method is customary and accepted in computability theory. Compare this proof with the proof of Theorem~4.3.11 in \cite{Wei00}. Via Lemma~4.3.6 it uses the closure of computable functions under primitive recursion (Theorem~3.1.7), which follows easily from Theorem~\ref{t4} in this article.

In mathematical practice, often a function such as addition on $\IQ$, is said to be computable ``by Church's Thesis''. In such a case, implicitly  fixed ``natural'', ``effective'' or ``canonical'' representations by finite or infinite strings are presupposed such that the function is computable w.r.t. these representations. Usually there is no disagreement about the meaning of ``natural'', ``effective''  or ``canonical''.  Often a canonical representation of a set $X$ can be defined up to equivalence by requiring  that
the functions and (or) predicates of the natural structure for $X$ must become computable,
and requiring additionally that the representation is maximal or minimal w.r.t. reducibility $\leq$ when indicated.

 Let us call  sets on which computability can be defined by canonical representations ``natural''. Examples of natural sets $Y$ are
 $\IN$, $\IB:=\IN^\IN$ (Baire space), $\Gamma^*$ and  $\Gamma^\IN$ (for finite $\Gamma$), $\IQ$, $\IQ^n$ and  $\bigcup_n\IQ^n$. Let us call a multi-representations $\delta :Y\mto X$
natural if $Y$ is natural.
By Theorem~\ref{t5},  Theorem~\ref{t6} can be generalized as follows.

\begin{thm}[informal generalization]\label{t6}
Let ${\bf P}$ be a generalized Turing machine on sets with natural representations.
Suppose that every function and test used in the machine has a  realization that is computable by Church's Thesis. Then the function $f_P$ computed by the machine is computable w.r.t. the natural representations.
\end{thm}

Brattka and Gherardi \cite{BG11} use a reduction $\leq_W$ for comparing the non-computability of theorems in analysis. Formally, $\leq_W$ compares the non-computability
of multi-functions on represented sets. We generalize this definition to multi-represented sets $(Z_i,\delta_i)$ ($i\in\{1,2,3,4\}$): For $f:Z_1\mto Z_2$ and $g:Z_3\mto Z_4$. $f\leq_W  g$, if there
are computable functions $G,H$ on $\om$ such that $p\mapsto G(p,h\circ H(p))$ is a
realization of $f$ if $h$ is a realization of $g$.
The next theorem provides a method to prove $f\leq_W g$.
We use the concept of {\em extension} for multi-functions from \cite[Definition~7]{Wei08}.
$f':Z_1\mto Z_2$ extends $f:Z_1\mto Z_2$, if
$\dom(f)\In\dom(f')$ and $f'(x)\In f(x)$ for all $x\in\dom(f)$.

For a generalized Turing machine $\bf M$ from Definition~\ref{d1} let
$\graph({\bf M}):=({\cal L}, S)$  where  $(l,l')\in S$ \ ($(l,l')$ is an edge), iff ${\rm Stm}(l)$ has the form $(\ldots, l')$ , $(\ldots {\rm then} \ l',\ {\rm else}\ l'')$ or
$(\ldots {\rm then} \ l'',\ {\rm else}\ l')$.

\begin{thm}\label{t9} For multi-functions $f,g$ on multi-represented sets such that $g$ is not computable, $f\leq_W  g$ \  if there is a generalized Turing machine $\bf N$ on represented sets such that
\begin{enumerate}[\em(1)]
\item \label{t9a} every test in $\bf N$ is computable,
\item \label{t9b} for every statement of the form $(i:=c(i_1,\ldots,i_n),l')$, either $c$ is computable or it is of the form $(i:=g(i_1),l')$.
\item \label{t9c}  every path in $\graph({\bf N})$ starting at $l_0$ visits at most once a label l such that ${\rm Stm}(l)$ applies the function $g$.
\item \label{t9d} $f_N$, the function computed by $\bf N$, extends $f$.
\end{enumerate}
(where ``computable'' means computable w.r.t the given multi-representations.)
\end{thm}

Condition (\ref{t9c}) for the GTM can be enforced easily syntactically.
The theorem generalizes  one direction of \cite[Lemma~4.5]{GM09}.

\pproof
Let $h\pf\om\to\om$ be a realization of $g$.
There is a generalized Turing machine $\bf M$ on $\om$ that realizes $\bf N$ (Definition~\ref{d5}) such that in $\bf M$ every test is computable, every function is computable or equal to $h$
and Condition (\ref{t9c}) is true for $\bf M$ and $h$.
By Theorem~\ref{t1}, $f_M$ realizes $f_N$ and hence $f_M$ realizes $f$.

For computing the functions $G$ and $H$, from the machine $\bf M$ we construct machines ${\bf M}_G$ and ${\bf M}_H$ .
Let ${\bf M}=({\mathcal L},l_0,l_f,\Gamma, b, 1,L,(X_i)_{0\leq i\leq L},{\rm Stm})$ where $X_i=\om$ for all $i$.

Define  ${\bf M}_H:=({\mathcal L},l_0,l_f,\Gamma, b, 1,L,(X_i)_{0\leq i\leq L},{\rm Stm}_H)$ such that for all $l,l',i,i_1$, %
\[{\rm Stm}_H(l):=\left\{\begin{array}{ll}
(0:={\rm id}_\om(i_1),l_f) &\mbox{if } \
{\rm Stm}(l)=(i:=h(i_1),l')\,,\\
{\rm Stm}(l) & \mbox{otherwise}\,.
\end{array}\right.\]
Then for input $p\in\dom(f_M)$, the machine ${\bf M}_H$ computes the argument for $h$ if a statement with $h$ is visited during the computation and computes the value $f_M(p)$ otherwise.

Define ${\bf M}_G:=({\mathcal L},l_0,l_f,\Gamma, b, 1,L+1,(X_i)_{0\leq i\leq {L+1}},{\rm Stm}_G)$ such that $X_{L+1}:=\om$ and for all $l,l',i,i_1$, %
\[{\rm Stm}_G(l):=\left\{\begin{array}{ll}
(i={\rm id}_\om(i_{L+1}),l')  &\mbox{if } \
{\rm Stm}(l)=(i:=h(i_1),l')\,,\\
{\rm Stm}(l) & \mbox{otherwise}\,.
\end{array}\right.\]
Then the machine ${\bf M}_G$ works in the same way as the machine ${\bf M}$ except for statements $(i:=h(i_1),l')$ of ${\bf M}$ where instead of applying $h$, ${\bf M}_G$ copies the value $q$ scanned by the head on Tape $L+1$ to the cell scanned by the head on Tape $i$.
For avoiding renaming of tapes we may assume w.l.o.g. that the machine ${\bf M}_G$ has the two  input tapes $1$ and $L+1$. Obviously for all $p\in\dom(f_M)$, $f_M(p)=f_{M_G}(p,h\circ f_{M_H}(p))$. Define $H:=f_{M_H}$ and $G:=f_{M_G}$.
\qq 
\medskip

Theorem~\ref{t9} holds accordingly for continuous reducibility instead of computable reducibility.

\section{Conclusion}\label{sech}
We have introduced the Generalized Turing machine as a simple general model of computation. This model is not intended for implementation on computers but as a mathematical tool for proving computability in Analysis. Although the three main
theorems \ref{t3}, \ref{t4} and \ref{t5} seem to be obvious and hence have already been applied informally without proofs, this article shows that even for our very meagre model of computation the proofs require some care.

\section*{Acknowledgement}
  The authors wish to thank the unknown referees for their careful work.

\bibliographystyle{plain}

\begin{thebibliography}{10}

\bibitem{Bla00}
Jens Blanck.
\newblock Domain representations of topological spaces.
\newblock {\em Theoretical Computer Science}, 247:229--255, 2000.

\bibitem{Blu04}
Lenore Blum.
\newblock Computing over the reals: where {T}uring meets {N}ewton.
\newblock {\em Notices of the AMS}, 51(9):1024--1034, 2004.

\bibitem{BCSS98}
Lenore Blum, Felipe Cucker, Michael Shub, and Steve Smale.
\newblock {\em Complexity and Real Computation}.
\newblock Springer, New York, 1998.

\bibitem{BSS89}
Lenore Blum, Mike Shub, and Steve Smale.
\newblock On a theory of computation and complexity over the real numbers:
  {$NP$}-completeness, recursive functions and universal machines.
\newblock {\em Bulletin of the American Mathematical Society}, 21(1):1--46,
  July 1989.

\bibitem{BC90}
Hans Boehm and Robert Cartwright.
\newblock Exact real arithmetic, formulating real numbers as functions.
\newblock In D.~Turner, editor, {\em Research topics in functional
  programming}, pages 43--64. Addison-Wesley, 1990.

\bibitem{Bra03f}
Vasco Brattka.
\newblock The emperor's new recursiveness: The epigraph of the exponential
  function in two models of computability.
\newblock In Masami Ito and Teruo Imaoka, editors, {\em Words, Languages \&
  Combinatorics III}, pages 63--72, Singapore, 2003. World Scientific
  Publishing.
\newblock ICWLC 2000, Kyoto, Japan, March 14--18, 2000.

\bibitem{BG11}
Vasco Brattka and Guido Gherardi.
\newblock Weihrauch degrees, omniscience principles and weak computability.
\newblock {\em Journal of Symbolic Logic}, 76:143--176, 2011.

\bibitem{BH98}
Vasco Brattka and Peter Hertling.
\newblock Feasible real random access machines.
\newblock {\em Journal of Complexity}, 14(4):490--526, 1998.

\bibitem{BHW08}
Vasco Brattka, Peter Hertling, and Klaus Weihrauch.
\newblock A tutorial on computable analysis.
\newblock In S.~Barry Cooper, Benedikt L\"owe, and Andrea Sorbi, editors, {\em
  New Computational Paradigms: Changing Conceptions of What is Computable},
  pages 425--491. Springer, New York, 2008.

\bibitem{BC06}
Mark Braverman and Stephen Cook.
\newblock Computing over the reals: Foundations for scientific computing.
\newblock {\em Notices of the AMS}, 53(3):318--329, 2006.

\bibitem{GM09}
Guido Gherardi and Alberto Marcone.
\newblock How incomputable is the separable {H}ahn-{B}anach theorem?
\newblock {\em Notre Dame Journal of Formal Logic}, 50(4):293--425, 2009.

\bibitem{Grz55}
Andrzej Grzegorczyk.
\newblock Computable functionals.
\newblock {\em Fundamenta Mathematicae}, 42:168--202, 1955.

\bibitem{Grz57}
Andrzej Grzegorczyk.
\newblock On the definitions of computable real continuous functions.
\newblock {\em Fundamenta Mathematicae}, 44:61--71, 1957.

\bibitem{HU79}
John~E. Hopcroft and Jeffrey~D. Ullman.
\newblock {\em Introduction to Automata Theory, Languages and Computation}.
\newblock Addison-Wesley, Reading, 1979.

\bibitem{KW85}
Christoph Kreitz and Klaus Weihrauch.
\newblock Theory of representations.
\newblock {\em Theoretical Computer Science}, 38:35--53, 1985.

\bibitem{Lac55e}
Daniel Lacombe.
\newblock Extension de la notion de fonction r\'{e}cursive aux fonctions d'une
  ou plusieurs variables r\'{e}elles {I-III}.
\newblock {\em Comptes Rendus Acad\'{e}mie des Sciences Paris},
  240,241:2478--2480,13--14,151--153, 1955.
\newblock Th\'{e}orie des fonctions.

\bibitem{Luc77}
Horst Luckhardt.
\newblock A fundamental effect in computations on real numbers.
\newblock {\em Theoretical Computer Science}, 5(3):321--324, 1977.

\bibitem{Sch03}
Matthias Schr\"oder.
\newblock Admissible representations for continuous computations.
\newblock Informatik Berichte 299, FernUniversit\"at Hagen, Hagen, April 2003.
\newblock Dissertation.

\bibitem{TZ04}
J.V. Tucker and J.I. Zucker.
\newblock Abstract versus concrete computation on metric partial algebras.
\newblock {\em ACM Transactions on Computational Logic}, 5(4):611--668, 2004.

\bibitem{Wei00}
Klaus Weihrauch.
\newblock {\em Computable Analysis}.
\newblock Springer, Berlin, 2000.

\bibitem{Wei08}
Klaus Weihrauch.
\newblock The computable multi-functions on multi-represented sets are closed
  under programming.
\newblock {\em Journal of Universal Computer Science}, 14(6):801--844, 2008.

\bibitem{ZW03}
Ning Zhong and Klaus Weihrauch.
\newblock Computability theory of generalized functions.
\newblock {\em Journal of the Association for Computing Machinery},
  50(4):469--505, 2003.

\end{thebibliography}

\vspace{-20 pt}
\end{document}